\documentclass[12pt,reqno]{amsart}
\usepackage{amssymb, amsfonts, amsbsy, latexsym, color}
\usepackage[dvips]{graphicx}

\newcommand{\bitem}{\begin{itemize}}
\newcommand{\eitem}{\end{itemize}}
\newcommand{\beq}{\begin{equation}}
\newcommand{\eeq}{\end{equation}}
\newcommand{\goto}{\rightarrow}
\newcommand{\cC}{{\mathcal C}}

\newcommand{\cA}{{\mathcal A}}
\newcommand{\cP}{{\mathcal P}}
\newcommand{\cS}{{\mathcal S}}

\newcommand{\bR}{{\bf R}}
\newcommand{\bP}{{\bf P}}
\newcommand{\bZ}{{\bf Z}}

\newcommand{\cR}{{\mathcal R}}
\newcommand{\cT}{{\mathcal T}}
\newcommand{\cH}{{\mathcal H}}
\newcommand{\cN}{{\mathcal N}}
\newcommand{\eps}{{\varepsilon}}

\newcommand{\SgSp}{\mbox{sing supp}}

\newcommand{\cL}{{\mathcal L}}
\newcommand{\OneStep}{{\sc One-Step}~}
\def\t{\tilde}
\newcommand{\ip}[2]{\langle#1,#2\rangle}

\newcommand{\absip}[2]{| \langle#1,#2\rangle |}
\newcommand{\norm}[1]{\|#1\|}

\definecolor{cardinal}{rgb}{.64,0.,.11}

\numberwithin{equation}{section}
\newtheorem{theorem}{Theorem}[section]
\newtheorem{corollary}{Corollary}[section]
\newtheorem{lemma}{Lemma}[section]
\newtheorem{proposition}{Proposition}[section]
\newtheorem{definition}{Definition}[section]

\begin{document}

\title[Geometric Separation by Single-Pass Alternating Thresholding]{Geometric Separation by\\ Single-Pass Alternating Thresholding}

\author[G. Kutyniok]{Gitta Kutyniok}
\address{Department of Mathematics, Technische Universit\"at Berlin,  10623 Berlin, Germany}
\email{kutyniok@math.tu-berlin.de}

\thanks{The author would like to thank David Donoho for numerous discussions on this and related topics.
She is grateful to the Department of Statistics at Stanford University and the Department
of Mathematics at Yale University for their hospitality and support during her visits, and
would also like to thank the Newton Institute of Mathematics in Cambridge, UK for
providing an inspiring research environment which led to the completion of a
significant part of this work during her stay. The author acknowledges support by the Einstein
Foundation Berlin, by Deutsche Forschungsgemeinschaft (DFG) Heisenberg fellowship KU 1446/8,
Grant SPP-1324 KU 1446/13 and DFG Grant KU 1446/14, and by the DFG Research Center {\sc Matheon}
``Mathematics for key technologies'' in Berlin.}

\begin{abstract}
Modern data is customarily of multimodal nature, and analysis tasks typically require separation into
the single components. Although a highly ill-posed problem, the morphological difference of these
components sometimes allow a very precise separation such as, for instance, in neurobiological
imaging a separation into spines (pointlike structures) and dendrites (curvilinear structures).
Recently, applied harmonic analysis introduced powerful methodologies to achieve this task,
exploiting specifically designed representation systems in which the components are sparsely
representable, combined with either performing $\ell_1$ minimization or thresholding on the combined dictionary.

In this paper we provide a thorough theoretical study of the separation of a distributional
model situation of point- and curvilinear singularities exploiting a surprisingly simple
single-pass alternating thresholding method applied to the two complementary frames:
wavelets and curvelets. Utilizing the fact that the coefficients are {\em clustered
geometrically}, thereby exhibiting {\em clustered/geometric sparsity} in the chosen
frames, we prove that at sufficiently fine scales arbitrarily precise separation is
possible. Even more surprising, it turns out that the thresholding index sets converge to
the wavefront sets of the point- and curvilinear singularities in phase space {\em and}
that those wavefront sets are perfectly separated by the thresholding procedure.
Main ingredients of our analysis are the novel notion of {\em cluster coherence} and
{\em clustered/geometric sparsity} as well as a {\em microlocal analysis viewpoint}.
\end{abstract}

\keywords{Thresholding. Sparse Representation.
Mutual Coherence. Tight Frames. Curvelets, Shearlets,
Radial Wavelets. Wavefront Set}

\maketitle

\section{Introduction}

Along with the deluge of data we face today, it is not surprising that the complexity of such
data is also increasing. One instance of this phenomenon is the occurrence of multiple
components, and hence, analyzing such data typically involves a separation step. One most
intriguing example comes from neurobiological imaging, where images of neurons from
Alzheimer infected brains are studied with the hope to detect specific artifacts of this
disease. The prominent parts of images of neurons are spines (pointlike structures)
and dendrites (curvelike structures), which require separate analyzes, for instance,
counting the number of spines of a particular shape, and determining the thickness of
dendrites \cite{KL12,MVVS04}.

From an educated viewpoint, it seems almost impossible to extract two images out of one image;
the only possible attack point being the morphological difference of the components. The new
paradigm of sparsity, which has lately led to some spectacular successes in solving such
underdetermined systems, does provide
a powerful means to explore this difference. The main sparsity-based approach towards solving such
separation problems consists in carefully selecting two representation systems, each one
providing a sparse representation of one of the components and both being incoherent with
respect to the other -- the encoding of the morphological difference --, followed by a procedure
which generates a sparse expansion in the dictionary combining the two representation systems.
This intuitively automatically forces the different components into the coefficients of the
`correct' representation system.

Browsing through the literature, the two main sparsity-based separation procedures can be
identified to be $\ell_1$ minimization (see, e.g.,
\cite{BGN08,DE03,DET06,DH01,DS89,EB02,ESQD05,GB03,GN03,KT09,SED04,SED05,SMBED05,ZP01})
and thresholding (see, e.g., \cite{BSFMD07,ESQD05,MZ93,MAC02}). For general papers on
$\ell_1$ minimization techniques we refer to \cite{CRT06b,CDS01,Don06a,Don06b,Don06c} and thresholding
to \cite{Tro04} or the reference list in the beautiful survey paper \cite{BDE09}.
While $\ell_1$ minimization has
produced very strong theoretical results, thresholding is typically significantly harder
to analyze due to its iterative nature. However, thresholding algorithms are in general much
faster than $\ell_1$ minimization, which makes them particularly attractive for the
aforementioned neurobiological imaging application due to its large problem size.

\medskip

In this paper we focus on thresholding as a separation technique for separating
point- from curvelike structures using radial wavelets and curvelets; in fact, we study
the very simple technique of single-pass alternating thresholding, which expands the image
in wavelets, thresholds and reconstructs the point part, then expands the residual in
curvelets, thresholds and reconstructs the curve part. In this paper we aim for a
fundamental mathematical understanding of the precision of separation allowed by this
thresholding method. Interestingly, our analysis
requires the notions of {\em cluster coherence} and {\em clustered/geometrical sparsity},
which were introduced in \cite{DK08a} by the author and Donoho in the context of analyzing $\ell_1$
minimization as a separation methodology.

We find the results in our paper quite surprising in two ways. First, the thresholding procedure
we consider is very simple, and researchers on thresholding algorithms might at first sight
dismiss such single-pass alternating thresholding methodology. Therefore, it is
intriguing to us, that we derive a quite similar perfect separation result (Theorem
\ref{theo:thresholding1}) as in our paper \cite{DK08a}, where $\ell_1$ minimization as a
separation technique was analyzed. Secondly, to our mind, it is even more surprising
that in Theorems \ref{theo:thresholding2} and \ref{theo:thresholding3} we derive even more
satisfying results by showing that the thresholding index sets converge to the wavefront sets
of the point- and curvilinear singularities in phase space {\em and} that those wavefront sets
are perfectly separated by the thresholding procedure. This, we already suspected for $\ell_1$
minimization to be true. However, we are not aware of any analysis tools strong enough to derive
these results for separation by $\ell_1$ minimization.


\subsection{A Geometric Separation Problem}

Let us start by defining the following simple but clear model problem of geometric separation  (compare also the
problem posted in \cite{DK08a}). Consider a `pointlike' object  $\cP$ made of point singularities:
\begin{equation} \label{pointdef}
   \cP = \sum_{i=1}^P |x - x_i|^{-3/2} .
\end{equation}
This object is smooth away from the $P$ given points $(x_i : 1 \leq i \leq P)$.
Consider as well a `curvelike' object
$\cC$, a singularity along a closed curve $\tau : [0,1] \mapsto \bR^2$:
\begin{equation} \label{curvedef}
    \cC = \int \delta_{\tau(t)}(\cdot)dt,
\end{equation}
where $ \delta_x$ is the usual Dirac delta function located  at $x$.
The singularities underlying these two distributions are geometrically quite different,
but the exponent $\-3/2$ is chosen so the  energy distribution across scales is similar;
if $\cA_r$ denotes the annular region $r < |\xi| < 2r$,
\[
      \int_{\cA_r} |\hat{\cP}|^2(\xi) \asymp r  ,
\qquad \int_{\cA_r} |\hat{\cC}|^2(\xi) \asymp r, \qquad r \goto \infty .
\]
This choice makes the components comparable
as we go to  finer scales;
the ratio of energies is more or less independent of scale.
 Separation is challenging at {\it every} scale.

Now assume that we observe
the `Signal'
\[
     f = \cP + \cC,
\]
however, the component distributions $\cP$ and $\cC$
are unknown to us.

\begin{definition}
The {\em Geometric Separation Problem}
requires to recover $\cP$ and $\cC$ from know\-ledge only of $f$;
here $\cP$ and $\cC$ are unknown to us, but obey (\ref{pointdef}),
(\ref{curvedef}) and certain regularity conditions on the curve $\tau$.
\end{definition}

As there are two unknowns ($\cP$ and $\cC$) and only
one observation ($f$), the problem seems improperly posed.
We develop a principled, rational approach which provably
solves the problem according to clearly stated standards.


\subsection{Two Geometric Frames }
\label{sec:twogeometricframes}

We now focus on two overcomplete systems  for representing the object $f$:
\bitem
\item {\it Radial Wavelets} -- a  tight frame with perfectly isotropic
generating elements.
\item {\it Curvelets} -- a highly directional tight frame with increasingly
anisotropic elements at fine scales.
\eitem
We pick these because, as is well known, point singularities
are coherent in the wavelet frame and curvilinear singularities
are coherent in the curvelet frame.
In Section \ref{sec:extensions} we discuss other system pairs. For readers
not familiar with frame theory, we refer to \cite{Chr03,CK12}, where terms
like `tight frame' --  a Parseval-like property --  are carefully discussed.

The point- and curvelike objects we defined in the previous subsection are real-valued distributions.
Hence, for deriving sparse expansions of those, we will consider radial wavelets and curvelets
consisting of real-valued functions. So  only angles associated with
radians $\theta \in [0,\pi)$ will be considered, which later on we will, as is customary, identify
with $\bP^1$, the real projective line.

We now construct the two selected tight frames as follows. Let $W(r)$ be an `appropriate' window function, where
in the following we assume that $W$ belongs to $C^\infty(\bR)$ and is compactly supported on
$[-2,-1/2] \cup [1/2,2]$ while being the Fourier transform of a wavelet. For instance, suitably scaled
Lemari\`{e}-Meyer wavelets possess these properties. We define {\em continuous radial wavelets}
at scale $a > 0$ and spatial position $b \in \bR^2$ by their Fourier transforms
\[
      \hat{\psi}_{a,b}(\xi) =  a \cdot W(a|\xi|) \cdot \exp\{ i b'\xi\}.
\]
The {\em wavelet tight frame} is then defined as a sampling of $b$ on a series of
regular lattices $\{ a_j  \bZ^2 \}$, $j \geq j_0$,
where $a_j = 2^{-j}$, i.e.,  the radial wavelets at scale $j$ and spatial position
$k = (k_1,k_2)'$ are given by the Fourier transform
\[
      \hat{\psi}_{\lambda}(\xi) =  2^{-j} \cdot W(|\xi|/2^{j}) \cdot \exp\{ i k'\xi /2^j\},
\]
where we let $\lambda = (j,k)$ index position and scale.

For the {\it same} window function $W$ and a `bump function' $V$, we define {\em continuous
curvelets} at scale $a>0$, orientation $\theta \in [0,\pi)$, and spatial position $b \in \bR^2$
by their Fourier transforms
\[
      \hat{\gamma}_{a,b,\theta}(\xi) =  a^{\frac{3}{4}}  \cdot W(a|\xi|) V(a^{-1/2}(\omega-\theta))
         \cdot \exp\{ i b' \xi \}.
\]
See \cite{CD04,CD05a} for more details.
The {\em curvelet tight frame} is then (essentially) defined as a sampling of $b$ on a
series of regular lattices
\[
   \{ R_{\theta_{j,\ell}} D_{a_j}  \bZ^2 \}, \qquad j \geq j_0, \quad \ell = 0, \dots, 2^{\lfloor j/2 \rfloor} -1 ,
\]
where $R_{\theta}$ is planar rotation by $\theta$ radians, $a_j = 2^{-j}$,
$\theta_{j,\ell} = \pi \ell / 2^{j/2}$, $\ell = 0, \dots, 2^{j/2}-1$,
and  $D_a$ is anisotropic dilation by $diag(a,\sqrt{a})$, i.e., the curvelets
at scale $j$, orientation $\ell$, and spatial position $k = (k_1,k_2)$ are
given by the Fourier transform
\[
      \hat{\gamma}_{\eta}(\xi) =  2^{-j\frac{3}{4}}  \cdot W(|\xi|/2^{j}) V((\omega-\theta_{j,\ell})2^{j/2})
         \cdot \exp\{ i (R_{\theta_{j,\ell}}D_{2^{-j}}k)' \xi \},
\]
where let $\eta = (j,k,\ell)$ index scale, orientation, and scale.
(For a precise statement, see \cite[Section 4.3, pp. 210-211]{CD05b}).

Roughly speaking, the radial wavelets are `radial bumps' with position $k/2^{j}$
and scale $2^{-j}$, while the curvelets live on anisotropic regions of width
$2^{-j}$ and length $2^{-j/2}$.  The wavelets are good at representing point singularities
while the curvelets are good at representing curvilinear singularities.

Using the {\it same} window $W$, we can construct a family of filters $F_j$ with transfer functions
\[
    \hat{F}_j(\xi) = W(|\xi|/2^{j}), \qquad \xi \in \bR^2 .
\]
These filters allow us to decompose a function $g$ into pieces $g_j$ with different scales,
the piece $g_j$ at subband $j$ arises from filtering $g$ using $F_j$:
\[
   g_j = F_j \star g;
\]
the Fourier transform $\hat{g}_j$ is supported in the annulus with inner radius $2^{j-1}$ and outer radius
$2^{j+1}$. Because of our assumption on $W$, we can reconstruct the original function from these
pieces using the formula
\[
   g = \sum_j F_j \star g_j, \qquad g \in L^2(\bR^2) .
\]

%
%
%

We now apply this filtering to our known image $f$, obtaining the truly geometric decomposition
\[
f_j = F_j \star f = F_j \star (\cP + \cC) = \cP_j + \cC_j
\]
for each scale $j$.

For future use,  let $\Lambda_j$ denote the collection of indices $(j,k)$ of wavelets
at level $j$, and
let $\Delta_j$ denote the indices $\eta = (j,k,\ell)$ of curvelets at level $j$.
%



\subsection{Separation via Thresholding}
\label{subsec:Theorem1}

We now consider a simple `one--step-thresholding' method -- which we also refer to as `single pass
alternating thresholding' method -- formalizing the first few steps of a recipe
for separation pointed out by Coifman and Wickerhauser \cite[Fig. 26(a-h)]{CW93} (cf. also \cite{DK08a}).
It is formally specified in Figure \ref{fig:onestepthresholdingCP}.

\begin{figure}[ht]
\centering
\framebox{
\begin{minipage}[h]{5.3in}
\vspace*{0.3cm}
{\sc \underline{One-Step-Thresholding}}

\vspace*{0.5cm}

{\bf Parameters:}\\[-2ex]
\begin{itemize}
\item Filtered signal $f_j$ for a scale $j$.\\[-2ex]
\item Thresholding parameter $\eps < 1/64$.
\end{itemize}

\vspace*{0.25cm}

{\bf Algorithm:}\\[-2ex]
\begin{itemize}
\item[1)] {\it Threshold Wavelet Coefficients:}\\[-2ex]
\begin{itemize}
\item[a)] Obtain wavelet coefficients $c_\lambda = \ip{f_j}{\psi_\lambda}$, $\lambda \in \Lambda_j$.\\[-2ex]
\item[b)] Apply threshold to obtain the set of significant coefficients $\cT_{1,j} = \{ \lambda : |c_\lambda| \ge 2^{\eps j}\}$.\\[-2ex]
\end{itemize}
\item[2)] {\it Reconstruct  Wavelet Component and Residualize:}\\[-2ex]
\begin{itemize}
\item[a)] Set $W_j = \sum_{\lambda \in \cT_{1,j}} c_\lambda \psi_\lambda$.\\[-2ex]
\item[b)] Set $\cR_j = f_j - W_j = \sum_{\lambda \in \cT_{1,j}^c} c_\lambda \psi_\lambda.$\\[-2ex]
\end{itemize}
\item[3)] {\it Threshold Curvelet Coefficients of Residual:}\\[-2ex]
\begin{itemize}
\item[a)] Compute $d_\eta = \ip{\cR_j}{\gamma_\eta}$ ,  $\gamma \in \Delta_j$.\\[-2ex]
\item[b)] Apply threshold to obtain the set of significant coefficients $\cT_{2,j} = \{ \eta : |d_\eta| \ge 2^{j(1/4-\eps)}\}.$\\[-2ex]
\end{itemize}
\item[4)] {\it Reconstruct Curvelet Component:}\\[-2ex]
\begin{itemize}
\item[a)] Compute $C_j = \sum_{\eta \in \cT_{2,j}} d_\eta \gamma_\eta.$
\end{itemize}
\end{itemize}

\vspace*{0.25cm}

{\bf Output:}\\[-2ex]
\begin{itemize}
\item Sets of significant coefficients: $\cT_{1,j}$ and $\cT_{2,j}$.\\[-2ex]
\item Approximations to $\cP_j$ and $\cC_j$: $W_j$ and $C_j$.
\end{itemize}
\vspace*{0.01cm}
\end{minipage}
}
\caption{{\OneStep} Thresholding Algorithm to approximately decompose $f_j = \cP_j + \cC_j$.}
\label{fig:onestepthresholdingCP}
\end{figure}

\OneStep is a very simple, easily implementable way to approximately decompose the signal $f_j$ into
purported pointlike and curvelike parts.
Currently popular thresholding algorithms are usually far more complex than \OneStep \hspace*{-0.15cm}:
they apply similar operations multiple times, with stopping rules, threshold adaptation, etc. It
therefore may be surprising that this very simple noniterative algorithm, with nonadaptive threshold, also
works well. The thresholds are even almost chosen as if the data wouldn't be composed at all:
The first threshold $2^{\eps j}$ is chosen coarsely below the decay rate $O(2^{j/2})$
of significant wavelet coefficients of the `naked' point singularity $\cP_j$;
the second threshold $2^{j(1/4-\eps)}$ is chosen just slightly below the decay rate $O(2^{j/4})$
of significant curvelet coefficients of the `naked' curvilinear singularity $\cC_j$. Notice
that we threshold the wavelet component more aggressively; and we refer to Section \ref{subsec:roadmap}
for more precise heuristics on the choice of these two thresholds.  It comes as a
second surprise that our estimates as well as the framework of geometric separation are
strong enough to survive this `brutally simple' thresholding strategy, as it is shown in
the following result as well as Theorems \ref{theo:thresholding2} and \ref{theo:thresholding3}.

For the following result, which will be proven in Section \ref{sec:thresholdingtheorem1},
we continue to suppose the sequence $(f_j)_j$ is known; thus the ideal decomposition
into a pointlike and curvelike part would be given by $f_j = \cP_j + \cC_j$. We apply
\OneStep\hspace*{-0.15cm}, which outputs approximations $W_j$ and $C_j$ to $\cP_j$ and
$\cC_j$, respectively.

\begin{theorem} {\sc Asymptotic Separation via \OneStep Thresholding.} \label{theo:thresholding1}
\[
\frac{ \| W_j - \cP_j \|_2 + \| C_j - \cC_j \|_2 }{\| \cP_j\|_2 + \|\cC_j\|_2 } \goto 0, \qquad j
\goto \infty.
\]
\end{theorem}

It is well-known that $\ell_1$ minimization and thresholding are closely connected in various ways.
In the past few years it has been frequently found that results on successful  $\ell_1$ minimization
subsequently inspired parallel results on thresholding methods. In a particular sense, this happened
here as well; after obtaining an asymptotic separation result using $\ell_1$ minimization (cf. \cite{DK08a}),
we found a similar result for this surprisingly simple thresholding procedure. However, even more intriguingly,
when performing this thresholding procedure -- as opposed to $\ell_1$ minimization -- we are able to even
derive much more satisfying results than Theorem \ref{theo:thresholding1}, which we turn our attention
to now.


\subsection{Wavefront Set Separation}
\label{subsec:Theorem2_3}

The very simplicity of \OneStep makes it possible to analyze
delicate phenomena which do not seem analytically
tractable for iterative thresholding or even for the $\ell_1$ minimization problem considered in \cite{DK08a}.

The geometric separation model we have been studying is
distinguished by  the behavior of its singularities.  One might hope
that the two purported geometric components
$\t{C}$ and $\t{P}$, defined by
\[
\t{P} = \sum_j F_j \star W_j \quad \mbox{and} \quad \t{C} = \sum_j F_j \star C_j,
\]
have exactly the singularities that one expects.  To articulate this goal
requires the notions of wavefront set and phase space
from microlocal analysis, which are reviewed  below
and in Section \ref{sec:microlocal}.  Intuitively, phase space is
the collection of location/direction pairs and the wavefront set $WF(f)$
of a distribution is the subset of phase space
where $f$ exhibits singularities.  Point singularities
are omnidirectional, while curvilinear singularities point in one direction.

Theorem \ref{theo:thresholding1} shows that the distributions
$\cP$ and $\cC$ can be arbitrarily well approximated by thresholding -- a similar
result was derived in our companion paper \cite{DK08a} for $\ell_1$ minimization.
However, the most desirable and also rhetorically effective matching
condition would be an arbitrarily perfect approximation also of the associated wavefront
sets $WF(\cP)$ and $WF(\cC)$.


Surprisingly, we derive two results in this direction for \OneStep -- one on the `analysis'
side and the other on the `synthesis' side. The first result shows that the wavefront sets of
$\cP$ and $\cC$ can indeed be approximated with arbitrary high precision by the significant thresholding
coefficients $\cT_{1,j}^{PS}$ and $\cT_{2,j}^{PS}$. As a measure of distance we employ the nonsymmetric
Hausdorff-style distance $d(A,B)$, say, in phase space measuring the largest distance from any point of
a subset of phase space $A$ to the closest corresponding point of a different subset $B$.
As a second result, we prove that the wavefront sets of the
synthesized objects $\sum_j F_j \star C_j$ and $\sum_j F_j \star P_j$ coincide with $WF(\cP)$ and $WF(\cC)$,
respectively. We might interpret this result as recovering $WF(\cP)$ and $WF(\cC)$
from the composed image $f$, hence in this sense we do not only separate the pointlike
structures from the curvelike structures, but even more separate their wavefront sets.

For a precise statement of the aforementioned two results, we require to introduce some notions
from microlocal analysis, which will be our main analysis methodology. Phase space is the space of
all direction/location pairs $(b,\theta)$, where $b \in \bR^2$ and the orientational component $\theta$
will be regarded as an element in $\bP^1$, the real projective space\footnote{Here we identify
$\bP^1$ with $[0,\pi)$ and freely write one or the other in what follows.
It may at first seem more natural to think of  directions $[0,2\pi)$ rather than
orientations $[0,\pi)$, note however that in this paper
we consider {\it real-valued} distributions  $\cP +\cC$ measured by real-valued
curvelets $\gamma_\eta$ so directions are not resolvable, only orientations.
We also frequently abuse notation as follows: we will write $|\theta - \theta'|$
when what is actually meant is geodesic distance between two points on $\bP^1$.} in $\bR^2$.

Since radial wavelets are oriented in all directions, we denote the set of significant phase space
pairs produced by the wavelet component of algorithm \OneStep by
\beq \label{eq:T1j_PS}
 \cT_{1,j}^{PS} = \{b_{j,k} : (j,k) \in \cT_{1,j}\} \times \bP^1;
\eeq
the set of significant phase space pairs for the curvelet component
of \OneStep is:
\beq \label{eq:T2j_PS}
\cT_{2,j}^{PS} = \{(b_{j,k,\ell},\theta_{j,\ell}) : (j,k,\ell) \in \cT_{2,j}\}.
\eeq
We further require the notion of a metric in phase space, which we choose to be
\[
d_{PS}((b,\theta),(b',\theta')) = \left(\norm{b-b'}_2^2+|\theta-\theta'|^2\right)^{1/2},\qquad (b,\theta),(b',\theta') \in \bR^2 \times \bP^1.
\]
and its associated asymmetric distance
\[
d_{PS}(C,C') = \max_{c \in C} \min_{c' \in C'} \norm{c-c'}_2,\qquad \mbox{where } C, C' \subseteq \bR^2 \times \bP^1.
\]

Section \ref{sec:approxWF} then proves the following theorem.

\begin{theorem} {\sc Approximation of the Wavefront Sets.} \label{theo:thresholding2}
\bitem
\item[{\rm (i)}]
\[
\limsup_{j \to \infty} \: d_{PS}(\cT_{1,j}^{PS},WF(\cP)) = 0.
\]
\item[{\rm (ii)}]
\[
\limsup_{j \to \infty} \: d_{PS}(\cT_{2,j}^{PS},WF(\cC)) = 0.
\]
\eitem
\end{theorem}

In short, the significant coefficients in each purported
geometric component  cluster increasingly around the wavefront set
of the underlying `true' geometric component. We further derive the
following result (proved in Section \ref{sec:sepWF}).

\begin{theorem} {\sc Separation of the Wavefront Sets.} \label{theo:thresholding3}
\[
WF(\sum_j F_j \star W_j) = WF(\cP)  \qquad \mbox{and} \qquad  WF(\sum_j F_j \star C_j) = WF(\cC).
\]
\end{theorem}

This implies that the wavefront sets of the reconstructed components
are precisely what we might hope for.

It seems plausible that results similar to Theorems \ref{theo:thresholding2} and
\ref{theo:thresholding3} could hold, in particular, also for separation via
$\ell_1$ minimization, but we don't know of analytical tools powerful enough to prove this.


\subsection{Extensions}
\label{sec:extensions}

We would like to point out that the analysis of \OneStep for solving the special
separation problem we focus on in this paper, gives rise to very extensive
generalizations and extensions; a few examples are stated in the sequel.

\bitem
\item {\it More General Classes of Objects.}
Theorems \ref{theo:thresholding1}--\ref{theo:thresholding2} can be generalized to other
situations.  First, we could consider singularities of different orders.
This would allow $\cC$ to model `cartoon' images, where
the curvilinear singularities are now the boundaries of the pieces
for piecewise $C^2$ functions. Second, we can allow smooth perturbations, i.e.,
$f = (\cP + \cC + g) \cdot h $ where $g, h$ are smooth functions of rapid decay at $\infty$.
In this situation, we let the denominator in Theorem \ref{theo:thresholding1} be simply $\| f_j \|_2$.

\item {\it Other Frame Pairs.}
Theorems \ref{theo:thresholding1}--\ref{theo:thresholding2} hold without change for many other
pairs of frames and bases, such as, e.g., by \cite{DK08b}, for the pair orthonormal separable Meyer wavelets
and shearlets (cf. \cite{GKL06,KL07,KKL10,KL10}).

\item {\it Noisy Data.}
Theorems \ref{theo:thresholding1}--\ref{theo:thresholding2} are resilient to noise impact; an image
composed of $\cP$ and $\cC$ with additive `sufficiently small' noise exhibits the same
asymptotic separation.

\item {\it Rate of Convergence.}
Theorem \ref{theo:thresholding1} can be accompanied by explicit decay
estimates.
\eitem


\section{Microlocal Analysis Viewpoint}
\label{sec:microlocal}

The morphological difference between the two structures we intend to extract
-- points and curve -- is the key to separation. In the section we will describe why heuristically this key
issue makes separation possible as well as present our main means to choose the `correct' thresholds.

\subsection{Point- and Curvelike Structures in Phase Space}

Our intuition as well as hard analysis is based on a microlocal analysis viewpoint, which through the notion of
wavefront sets will allow us to, roughly speaking, include the morphology of the structures by adding a third
dimension to spatial domain. Let us  start by recalling the notion of wavefront set and  -- related with this
-- the notion of singular supports and phase space. The {\em singular support} of a distribution $f$, $\SgSp(f)$, is
defined to be the set of points where $f$ is not locally $C^\infty$. The notion of wavefront set then goes beyond
the classical spatial domain picture and extends it to {\em phase space}, which consists of position-orientation
pairs $(b,\theta)$; see the more detailed discussion in Section \ref{subsec:Theorem2_3}.
The {\em wavefront set} $WF(f)$ lives in this phase space and can be coarsely described
as the set of position-orientation pairs at which $f$ is nonsmooth; for more details, see: \cite{Hoe03,CD05a,KL07}.

\begin{figure}[ht]
\centering
\includegraphics[height=2.3in]{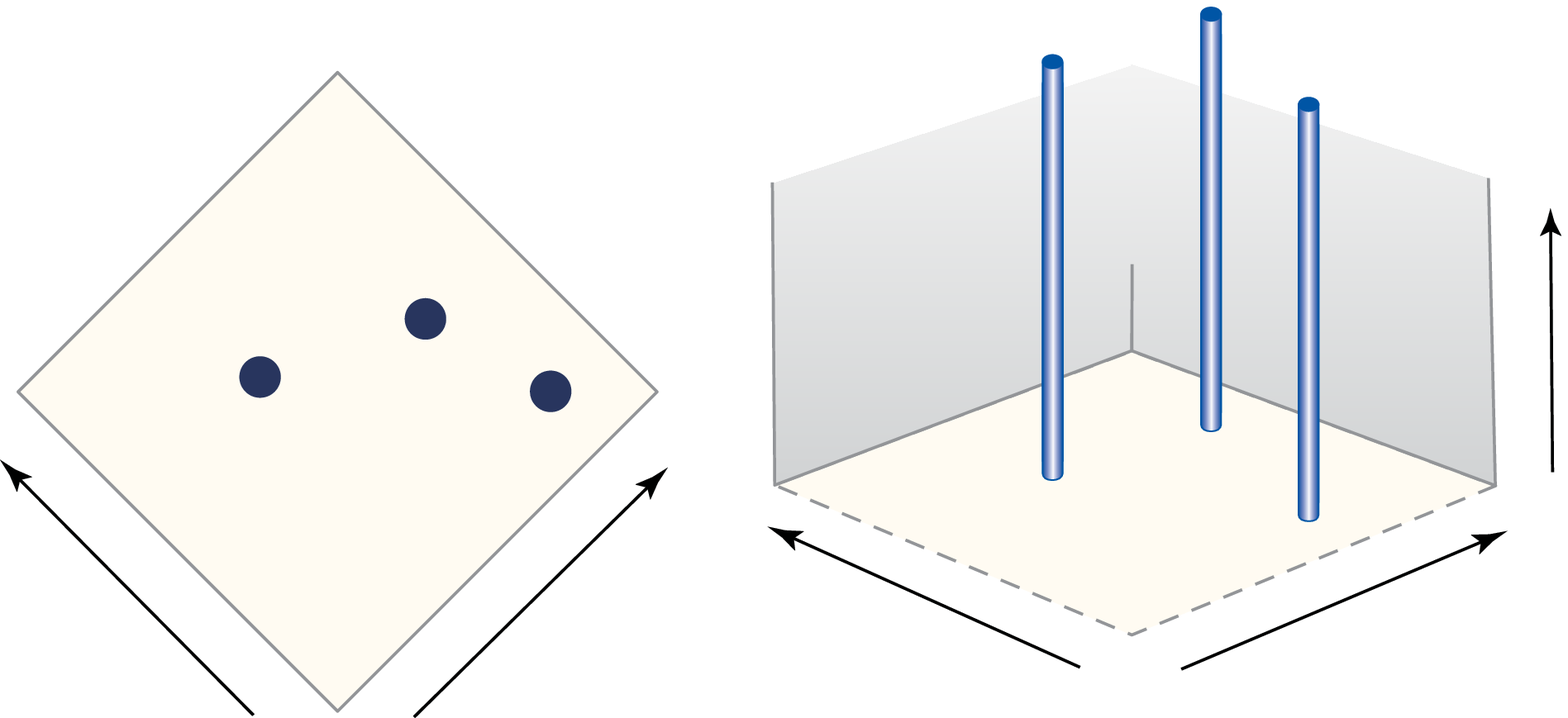}
\put(-367,45){$b_2$}
\put(-208,45){$b_1$}
\put(-185,30){$b_2$}
\put(-15,30){$b_1$}
\put(3,110){$\theta$} \caption{Left panel: singular support of
$\cP$. Right panel: wavefront set of  $\cP$ in phase
space.} \label{fig:PinPhaseSpace}
\end{figure}

To illustrate these notions and also prepare our heuristic argument why separation through thresholding is
possible, we first consider the distribution $\cP$. A short computation shows that
\[
\SgSp(\cP) = \{ x_i \} \qquad \mbox{and} \qquad WF(\cP) = \SgSp(\cP) \times \bP^1,
\]
which can be regarded as a manifestation of the isotropic nature of the point singularities. Illustrations of
$\SgSp(\cP)$ and of $WF(\cP)$ are presented in Figure \ref{fig:PinPhaseSpace}.

\begin{figure}[ht]
\centering
\includegraphics[height=2.3in]{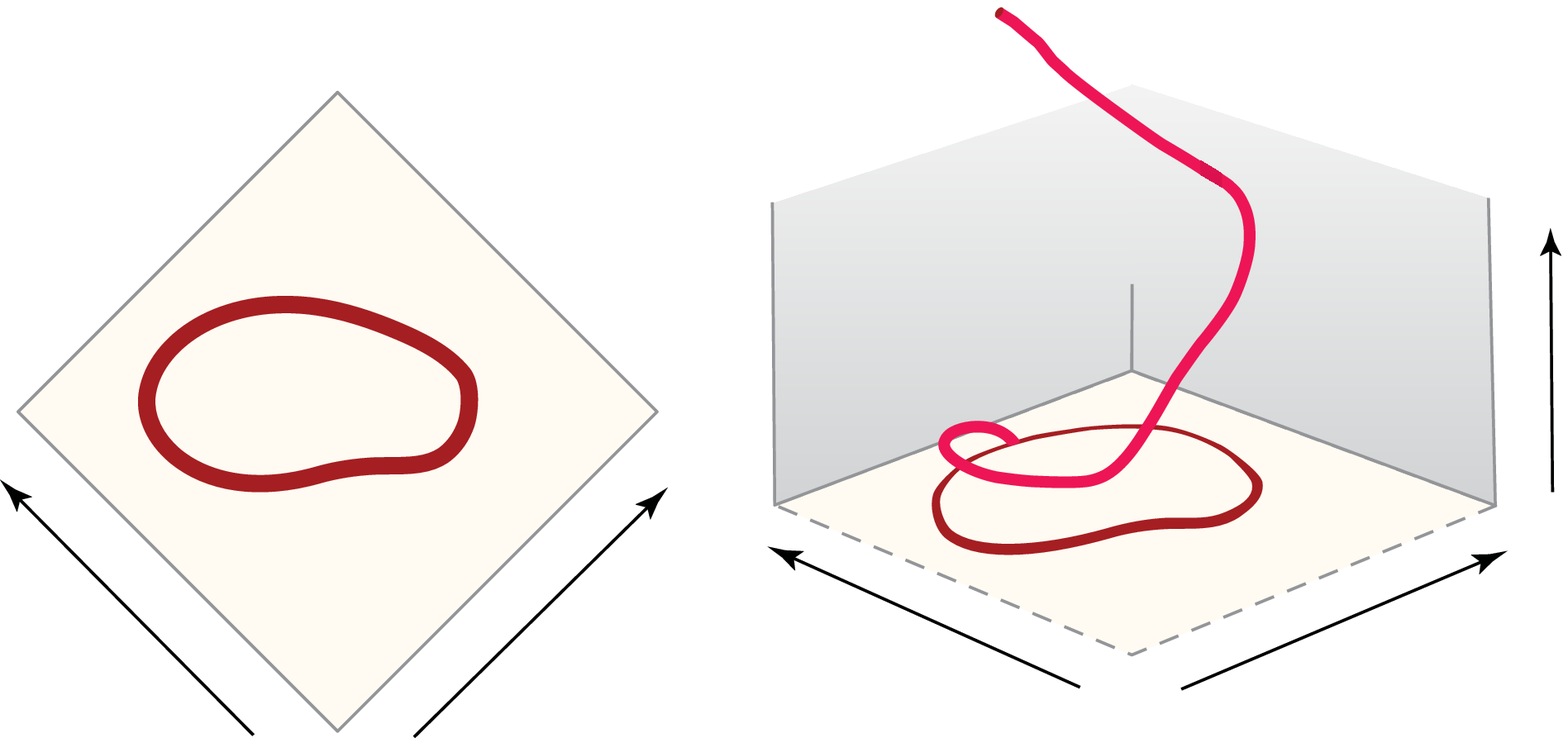}
\put(-360,45){$b_2$}
\put(-203,45){$b_1$}
\put(-180,30){$b_2$}
\put(-15,30){$b_1$}
\put(3,110){$\theta$} \caption{Left panel: singular support of
$\cC$. Right panel: wavefront set of  $\cC$ in phase
space.} \label{fig:CinPhaseSpace}
\end{figure}

For the distribution $\cC$, we obtain
\[
  \SgSp(\cC) = image(\tau) \qquad \mbox{and} \qquad WF(\cC) =  \{ (\tau(t), \theta(t)): t \in [0,L(\tau)] \},
\]
where $\tau(t)$ is a unit-speed parametrization of $\cC$ and $\theta(t)$ is the normal direction to $\cC$ at
$\tau(t)$ regarded in $\bP^1$. Here, the anisotropy and -- in comparison with Figure \ref{fig:PinPhaseSpace} --
the morphological difference to $\cP$ becomes evident. An illustration of $\SgSp(\cC)$ and of $WF(\cC)$
is presented in Figure \ref{fig:CinPhaseSpace}.

\subsection{Wavelets and Curvelets in Phase Space}

Although being smooth functions, in a certain sense, wavelets and curvelets can be regarded as leaving an
approximate footprint in phase space. To make this statement rigorous, we first observe the approximate footprint
in spatial domain left by wavelet and curvelets as detailed in the following two lemmata taken from \cite{DK08a}.
As expected, these observations show the isotropic nature of wavelets in contrast to the anisotropic nature
of curvelets.

\begin{lemma}[\cite{DK08a}]
\label{lemma:estimate_psi}
 For each $N = 1,2, \dots$
there is a constant $c_N$ so that
\[
   |\psi_{a,b}(x)| \leq c_N \cdot  a^{-1} \cdot  \langle | x - b |/a \rangle^{-N}, \qquad \forall a \in \bR^{+} \;
   \forall   b,x \in \bR^2.
\]
\end{lemma}

\begin{lemma}[\cite{DK08a}]
\label{lemma:estimate_gamma}
For each $N = 1,2, \dots$
there is a constant $c_N$ so that
\begin{equation*}
   |\gamma_{a,b,\theta}(x)| \leq c_N \cdot  a^{-3/4} \cdot  \langle | x - b |_{a,\theta} \rangle^{-N},
   \qquad \forall a \in \bR^{+} \;  \forall \theta \in [0,\pi) \;
   \forall   b,x \in \bR^2.
\end{equation*}
\end{lemma}

\begin{figure}[ht]
\centering
\includegraphics[height=2.3in]{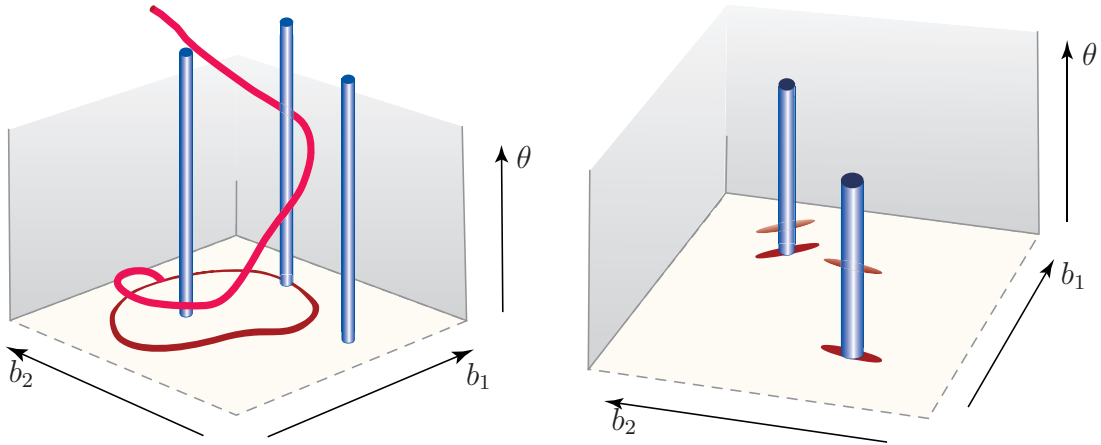}
\put(-403,23){$b_2$}
\put(-230,22){$b_1$}
\put(-211,104){$\theta$}
\put(-175,5){$b_2$}
\put(-5,60){$b_1$}
\put(3,143){$\theta$} \caption{Left panel: wavefront set of the observed data $f=\cP + \cC$.
Right panel: phase space footprint of radial wavelets and curvelets.} \label{fig:framesinPhaseSpace}
\end{figure}

Since it is known from \cite{CD05a} that the continuous curvelet transform precisely resolves the wavefront set of
distributions, we might consider the image of wavelets and curvelets under the continuous curvelet transform
for `sufficiently small' scale as a footprint of these in phase space. An illustration is given in Figure
\ref{fig:framesinPhaseSpace}, and for a detailed description we refer the interested reader to \cite{DK08a}.

Visually, wavelets are perfectly adapted to strongly react to $\cP$ in a similar way as curvelets
will strongly react to $\cC$. This will be now made precise and will lead to the chosen thresholds
for separation.

\subsection{Road Map to the `Correct' Thresholds}
\label{subsec:roadmap}

To slowly approach a rigorous phrasing of the aforementioned strong reaction, we first consider
the reaction of both wavelets to $\cP$ and $\cC$. A simplified form of Lemma \ref{lemma:pointwavelet}
states that,
\beq \label{eq:intuition1}
   \absip{\psi_{a_{j},x_i}}{\cP_j} = O(2^{j/2}) \qquad j \to \infty,
\eeq
with fast decay\footnote{As it is custom, we refer to the behavior $O(a^N)$ as $a \to 0$ for all $N = 1, 2, ...$ as {\em fast decay}.}
for other locations than $x_i$, and Lemma \ref{lemma:curvewavelet} shows that, for each
$j$ and $b$,
\beq \label{eq:intuition2}
|\langle \psi_{a_j,b}, \cC_j \rangle | =O(1).
\eeq
Secondly, turning our attention to curvelets and their reaction to $\cP$ and $\cC$, we observe that,
by a simplified form of Lemma \ref{lemma:pointcurvelet}, for each $\theta$,
\beq \label{eq:intuition3}
   \absip{\gamma_{a_{j},x_i,\theta}}{\cP_j} = O(2^{j/2})  \qquad j \to \infty,
\eeq
and, for $b$ positioned on the curve $\tau$ and $\theta$ pointing in the direction perpendicular to the tangent
to the curve in $b$,
\beq \label{eq:intuition4}
   \absip{\gamma_{a_{j},b,\theta}}{\cC_j} = O(2^{j/4}) \qquad j \to \infty.
\eeq

Examining closely \eqref{eq:intuition1} and \eqref{eq:intuition2}, it becomes immediately evident that
the correct first threshold -- which should capture $\cP_j$ by thresholding wavelet coefficients of
$f_j$ -- need to be chosen `slightly higher' than a constant asymptotically, wherefore
we choose it equal to $2^{\eps j}$ for small $\eps$.

The second threshold seem to be a somehow more serious problem, since \eqref{eq:intuition3} and \eqref{eq:intuition4}
show that curvelets react stronger to a point singularity than a curvilinear singularity. However, we wish
the reader to keep in mind that ideally all energy from $\cP_j$ is already captured during the first thresholding
procedure. Hence, it should presumably be `safe' to choose the second threshold -- which shall capture $\cC_j$ by thresholding
curvelet coefficients of the residual $\cR_j$ -- with asymptotic behavior $o(2^{j/4})$. To avoid unnecessary
risks, we choose it only slightly below $2^{j/4}$, more precisely, equal to $2^{j(1/4-\eps)}$.

\subsection{What Type of Separation Result is Preferable?}

Applying now \OneStep (cf. Figure \ref{fig:onestepthresholdingCP}) yields significant coefficient sets $\cT_{1,j}$ and
$\cT_{2,j}$ and approximations to $\cP_j$ and $\cC_j$: $W_j$ and $C_j$. In Sections \ref{subsec:Theorem1} and
\ref{subsec:Theorem2_3}, we presented three theorems on the `quality' of this separation, which we would
now like to discuss and compare.

Theorem \ref{theo:thresholding1} studies the relative separation error and proves that asymptotically this error
can be made arbitrarily small for sufficiently fine scale. This is in a sense the most natural question to
ask, and the theorem provides the answer one would hope for.

However, from a microlocal analysis viewpoint, the most satisfying separation to derive would be the perfect separation
of the wavefront sets of $\cP$ and $\cC$, i.e., to separate the LHS of Figure \ref{fig:framesinPhaseSpace} into
the RHS of Figures \ref{fig:PinPhaseSpace} and \ref{fig:CinPhaseSpace}. This would be considerably `stronger' than
Theorem \ref{theo:thresholding1} in the following sense: Once the wavefront sets are extracted, we have
complete information about the underlying singularities, in contrast to the merely asymptotic knowledge provided
by Theorem \ref{theo:thresholding1}.

Knowledge about $WF(\cP)$ and $WF(\cC)$ could be either coming from
$\cT_{1,j}$ and $\cT_{2,j}$ or from $W_j$ and $C_j$. The sets of significant coefficients generated by
thresholding do not provide an immediate means for separating the wavefront sets, since they live on the
analysis side (as opposed to the synthesis side). Astonishingly, they are still able to precisely {\em locate}
the wavefront sets $WF(\cP)$ and $WF(\cC)$, more precisely, they `converge' to the wavefront sets in phase
space measured in the phase space norm $d_{PS}$ as $j \to \infty$, which is the statement of Theorem
\ref{theo:thresholding2}. This shows that the points in $\cT_{1,j}$ and $\cT_{2,j}$ are located in
tubes around $WF(\cP)$ and $WF(\cC)$, respectively, which become more concentrated around these wavefront
sets as the scale becomes finer. The
corresponding objects on the synthesis side, i.e., $W_j$ and $C_j$, now allow separation of $WF(\cP)$ and
$WF(\cC)$, in the sense that the wavefront sets of the reconstructed distributions $\sum_j F_j \star W_j$
and $\sum_j F_j \star C_j$ precisely coincide with $WF(\cP)$ and $WF(\cC)$. This is the content of
Theorem \ref{theo:thresholding3}.


\section{Geometry of the Thresholded Wavelet Coefficients}
\label{sec:geometry_wavelets}

Following the ordering of the thresholding, we first focus on the set of significant radial wavelet
coefficients $\cT_{1,j}$ generated by Step 1) of {\sc One-Step-Thresholding} (see Figure \ref{fig:onestepthresholdingCP}),
in particular, on its phase space footprint, defined in \eqref{eq:T1j_PS} as
\[
\cT_{1,j}^{PS} = \{b \in \bR^2  : \absip{f_j}{\psi_{a_j,b}} \ge a_j^{-\eps}\}\times \bP^1.
\]
Our objective will be to derive a tube around $\cT_{1,j}^{PS}$ in phase space with controllable `size'. This
tube should therefore be a neighborhood of $WF(\cP)$, and hence be isotropic.

For our analysis, we first notice that WLOG we can assume that
\beq \label{eq:WLOG_P}
\cP=|x|^{-3/2}.
\eeq
From here, the result for the original $\cP$ as defined in \eqref{pointdef} can be concluded because
of the following reasons: Firstly, all results are translation invariant, hence instead of the origin
the results follow immediately for a different point in spatial domain; and secondly, the change from
one point to finitely many points just introduces a constant independent on $j$.


\subsection{Estimates for Wavelet Coefficients}

We start by analyzing the interaction of wavelet atoms. The technical proof of the  following result is provided
in Section \ref{subsec:proofs_microlocal}

\begin{lemma} \label{lemm:psi_psi_estimate}
For each $N=1,2,\dots$,
there is a constant $c_N$ so that
\[
|\langle \psi_{a,b} , \psi_{a_0,b_0} \rangle | \leq c_N \cdot 1_{\{|\log_2(a/a_0)| < 3\}} \cdot
    \langle |b - b_0|/a \rangle^{-N} .
\]
\end{lemma}

Next, we recall a result derived in \cite{DK08a} for radial wavelet coefficients of our
point singularity \eqref{eq:WLOG_P}.

\begin{lemma}[\cite{DK08a}]
\label{lemma:pointwavelet}
For each $N = 1,2, \dots$, there is a constant $c_N$ so that
\[
   \absip{\psi_{a_{j},b}}{\cP_j} \le c_N \cdot a_j^{-1/2} \cdot \langle |b/a_j| \rangle^{-N},
   \qquad \forall j \in \bZ\; \forall b \in \bR^{2}.
\]
\end{lemma}

In the sequel, we will further require an estimate of the wavelet coefficients of the curvilinear singularity $\cC_j$.
Notice that the following estimate does only provide a very coarse upper bound. In order to derive
a more detailed estimate, the curve would need to be much more carefully analyzed as it will be
done in Section \ref{sec:geometry_curvelets}. However, the estimate as stated below is all we will
require.


\begin{lemma}
\label{lemma:curvewavelet}
There exists a constant $c$ so that
\[
|\langle \psi_{a_j,b}, \cC     \rangle | \le c, \qquad \forall j \in \bZ\; \forall b \in \bR^{2}.
\]
\end{lemma}

\noindent
{\bf Proof.}
By Lemma \ref{lemma:estimate_psi} and the definition of the distribution $\cC$,
\beq \label{eq:psiC1}
|\langle \psi_{a_j,b}, \cC     \rangle |
\le \int_0^1 |\psi_{a_j,b}(\tau(t))|\,  dt
\le c_N \cdot  a^{-1} \cdot  \int_0^1  \langle | \tau(t) - b |/a \rangle^{-N} dt.
\eeq
WLOG we assume that $b \in \tau([0,1])$ with $\tau(0)=b$, say, and we can also assume that $b=(b_1,0)$. Choosing a ball $B_r(b)$ around
$b$ with $r$ chosen arbitrarily small (yet, independent of $j$), there exists some $0 < \delta < 1/2$ such that
\[
\tau([0,1]) \cap B_r(b) = \tau([1-\delta,1] \cup [0,\delta]).
\]
This information is now used to split the last integral in \eqref{eq:psiC1} according to
\beq \label{eq:psiC2}
\int_0^1  \langle | \tau(t) - b |/a \rangle^{-N} dt
= \int_{[1-\delta,1] \cup [0,\delta]}  \langle | \tau(t) - b |/a \rangle^{-N} dt
+ \int_{[\delta,1-\delta]}  \langle | \tau(t) - b |/a \rangle^{-N} dt
=: I_1 + I_2.
\eeq
For estimating $I_1$, we first observe that it is sufficient to consider $\int_{[0,\delta]}$ due to
symmetry reasons. For $r$ small enough, the curve inside $B_r(b)$ can be arbitrarily well approximated by its
osculating circle with its center denoted by $z=(z_1,0)$. Combining these considerations as well as exploiting
the approximation by a Taylor series for cosine,
\begin{eqnarray} \nonumber
\int_{[0,\delta]} \langle | \tau(t) - b |/a \rangle^{-N} dt
& \le & c \cdot \int_{[0,\delta]} \langle |(b_1-z_1)(\cos(t),\sin(t))+z - b |/a \rangle^{-N} dt\\ \nonumber
& = & c \cdot \int_{[0,\delta]} \langle |b_1-z_1| \sqrt{2(1-\cos(t))}/a \rangle^{-N} dt\\ \nonumber
& \le & c' \cdot \int_{[0,\delta]} \langle |b_1-z_1| \cdot t/a  \rangle^{-N} dt\\ \nonumber
& = & c' \cdot a/|b_1-z_1| \int_{[0,|b_1-z_1|\delta/a]} \langle t \rangle^{-N} dt\\ \label{eq:psiC3}
& \le & c'' \cdot a/|b_1-z_1|.
\end{eqnarray}
Using the definition of $r$, the integral $I_2$ can be easily estimated as
\beq \label{eq:psiC4}
\int_{[\delta,1-\delta]} \langle | \tau(t) - b |/a \rangle^{-N} dt
 \le  \int_{[\delta,1-\delta]} \langle r/a \rangle^{-N} dt
 \le  (r/a)^{-N}.
\eeq
Summarizing, by \eqref{eq:psiC1}--\eqref{eq:psiC4}, there exists some constant $c$ (independent on $a$ and $b$)
such that
\[
|\langle \psi_{a_j,b}, \cC     \rangle | \le c_N \cdot  a^{-1} \cdot (a/|b_1-z_1|+(r/a)^{-N}) \le c.
\]
The lemma is proved.
\qed


\subsection{Geometry of $\cT_{1,j}^{PS}$}

We now first analyze the set $\cT_{1,j}^{PS}$ by the following two lemmata.

\begin{lemma}
\label{lemma:T1c}
Let $(b,\theta) \in (\cT_{1,j}^{PS})^c$, and let $j$ be sufficiently large. Then, for each $N = 1,2, \dots$,
there is a constant $c_N$ so that
\[
|b/a_j| > c_N \cdot 2^{j\frac{1-2\eps}{2N}}.
\]
\end{lemma}

\noindent
{\bf Proof.}
Let $b \in \bR^2$ be such that
\[
|\ip{\psi_{a_j,b}}{\cP_j} + \ip{\psi_{a_j,b}}{\cC_j}| < 2^{j \eps}.
\]
Since by Lemma \ref{lemma:curvewavelet}, $|\ip{\psi_{a_j,b}}{\cC}|$, and hence $|\ip{\psi_{a_j,b}}{\cC_j}|$, is bounded by a
constant $c$, say, we have
\[
|\ip{\psi_{a_j,b}}{\cP_j}| < 2^{j \eps}+c.
\]
Next we use the estimate in Lemma \ref{lemma:pointwavelet} as a model to conclude that
\[
\langle |b/a_j| \rangle^{-N} < c_N \cdot 2^{-j/2} (2^{j \eps}+c).
\]
Thus, since for sufficiently large $j$, we have $2^{j \eps} > c$,
\[
|b/a_j| > \left((c_N \cdot 2^{-j/2} (2^{j \eps}+c))^{-2/N}-1\right)^{1/2} > \left(c_N \cdot 2^{j\frac{1-2\eps}{N}}-1\right)^{1/2}.
\]
Letting $j$ be large enough so that $(c_N/2) \cdot 2^{j\frac{1-2\eps}{N}} > 1$ proves the lemma.
\qed

\begin{lemma}
\label{lemma:T1}
Let $(b,\theta) \in \cT_{1,j}^{PS}$. Then, for each $N = 1,2, \dots$, there is a
constant $c_N$ so that
\[
|b/a_j| \le c_N \cdot 2^{j\frac{1-2\eps}{2N}}.
\]
\end{lemma}

\noindent
{\bf Proof.}
Let $b \in \bR^2$ be such that
\[
|\ip{\psi_{a_j,b}}{\cP_j} + \ip{\psi_{a_j,b}}{\cC_j}| \ge 2^{j \eps}.
\]
Since by Lemma \ref{lemma:curvewavelet}, $|\ip{\psi_{a_j,b}}{\cC}|$, and hence $|\ip{\psi_{a_j,b}}{\cC_j}|$,
is bounded by a constant $c$, say, we have
\[
|\ip{\psi_{a_j,b}}{\cP_j}| \ge 2^{j \eps}-c.
\]
Next we use the estimate in Lemma \ref{lemma:pointwavelet} as a model to conclude that
\[
\langle |b/a_j| \rangle^{-N} \ge c_N \cdot 2^{-j/2} (2^{j \eps}-c).
\]
Thus, since for sufficiently large $j$, we have $2^{j \eps-1} > c$,
\[
|b/a_j| \le \left((c_N \cdot 2^{-j/2} (2^{j \eps}-c))^{-2/N}-1\right)^{1/2} \le \left(c_N \cdot 2^{j\frac{1-2\eps}{N}}\right)^{1/2}.
\]
The lemma is proved.
\qed

We certainly hope (and expect) that the threshold is set in such a way that the wavefront set of $\cP$ is
contained in $\cT_{1,j}^{PS}$. This is obviously the first requirement for being able to separate both wavefront
sets $WF(\cP)$ and $WF(\cC)$ through Single-Pass Alternating Thresholding (compare Theorem \ref{theo:thresholding3}).
The next result shows that this is indeed the case.

\begin{lemma} \label{lemm:WFPinT1}
For $j$ sufficiently large,
\[
\absip{\cP_j}{\psi_{a_j,0}} \ge c \cdot 2^{j/2}.
\]
Hence, in particular,
\[
WF(\cP) \subseteq \cT_{1,j}^{PS}.
\]
\end{lemma}

\noindent
{\bf Proof.}
By Parseval,
\[
\absip{\cP_j}{\psi_{a_j,0}}
= 2\pi \absip{\hat{\cP}_j}{\hat{\psi}_{a_j,0}}
= 2\pi \cdot 2^{j/2} \int W^2(|\xi|) |\xi|^{-1/2} d\xi.
\]
Hence, we can conclude that, for $j$ sufficiently large,
\[
\absip{\cP_j}{\psi_{a_j,0}} \ge c \cdot 2^{j/2}.
\]
This proves the first claim.

For the `in particular'-part, recall that $WF(\cP) = \{0\} \times \mathbf{P}^1$.
By Lemma \ref{lemma:curvewavelet} and the previous consideration, for
sufficiently large $j$,
\[
\absip{f_j}{\psi_{a_j,0}}
\ge \absip{\cP_j}{\psi_{a_j,0}} - \absip{\cC_j}{\psi_{a_j,0}}
\ge \absip{\cP_j}{\psi_{a_j,0}} - c
\ge 2^{\eps j}.
\]
The lemma is proved.
\qed


\section{Geometry of the Thresholded Curvelet Coefficients}
\label{sec:geometry_curvelets}

This section now aims to derive a fundamental geometric understanding of the cluster of curvelet
coefficients $\cT_{2,j}$ generated by Step 3) of {\sc One-Step-Thresholding} of the residual generated in
Step 2) (see Figure \ref{fig:onestepthresholdingCP}).
The phase space geometry will play an essential role in setting up the analysis correctly, hence
it will be beneficial to study the projection of $\cT_{2,j}$ onto phase space, defined in
\eqref{eq:T2j_PS}, as
\[
\cT_{2,j}^{PS} = \{(b,\theta) \in \bR^2 \times [0,\pi) : \absip{\cR_j}{\gamma_{a_j,b,\theta}} \ge a_j^{\eps-1/4}\}.
\]
Morally, the points in phase space associated with significant curvelet
coefficients, given by $\cT_{2,j}^{PS}$, are contained in a tube around $WF(w\cL)$ in phase space.
The main objective will now be to explicitly define such a tube around the phase space footprint of this cluster, where we
have more control on. This will become crucial for handling the thresholded curvelet coefficients
in the proofs of Theorems \ref{theo:thresholding1}--\ref{theo:thresholding3}.




\subsection{Bending the Curve}
\label{subsec:bending}

We first face the problem of how to deal with the curvilinear singularity. In \cite{DK08a}, this problem was tackled by
carefully and smoothly breaking the curve into pieces, bending each piece, and then combining pieces in the end.
This technique shall also be applied here. For the convenience of the reader, we review the main ideas of this
particular approach.

First, a quantitative `tubular neighborhood theorem' is being developed to allow local bending of the curve. Due to
regularity of the curve, there exists some $\rho$ small compared to the curvature of $\tau$, so that
\[
     \int_{(i-1)\rho}^{(i+1)\rho} |\tau^{''} (t) | dt \leq \eps, \qquad i=0, ..., \tfrac{{\rm length}(\tau)}{\rho} =: m.
\]
Now consider the following local coordinate system in the vicinity of $\tau$. Let $t_i =  i \rho$, for $i=0 , \dots , m$
with $\tau(t_0) = \tau( t_m)$, since $\tau$ is closed. Then we have the following

\begin{figure}
\centering
\includegraphics[height=1.75in]{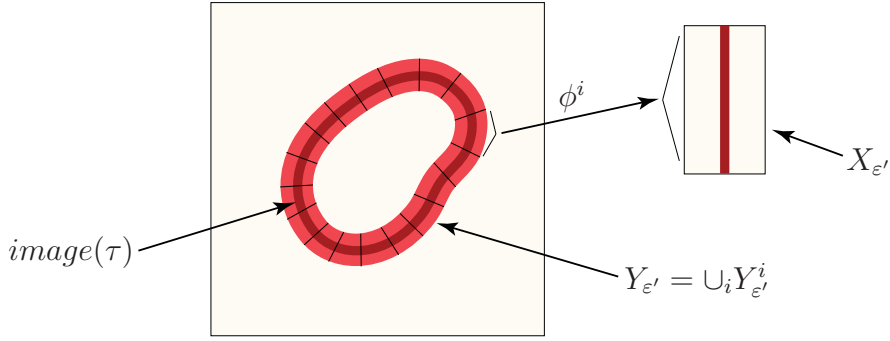}
\put(-316,29){$image(\tau)$}
\put(-83,18){$Y_{\eps'} = \cup_i Y_{\eps'}^i$}
\put(-108,87){$\phi^i$}
\put(1,63){$X_{\eps'}$}
\caption{The tubular neighborhood $Y_{\eps'} = \cup_i Y_{\eps'}^i$ of $image(\tau)$ and the
mapping $\phi^i : Y_{\eps'}^i \mapsto X_{\eps'}$.}
\label{fig:tube}
\end{figure}

\begin{lemma}[\cite{DK08a}] {\bf (Tubular Neighborhood Theorem)}
For sufficiently small $\eps > 0$, there is some $\eps' > 0$ so that, for
$X_{\eps'} = [-\eps',\eps'] \times [-\rho,\rho]$, we have:
\bitem
\item  for each $i=0 , \dots , m$, there exists a tube $Y_{\eps'}^i$ around $\tau$ and
an associated diffeomorphism $\phi^i : Y_{\eps'}^i \mapsto X_{\eps'}$,
\item  the mapping $\phi^i $ extends to a diffeomorphism from $\bR^2$ to $\bR^2$ which reduces to the identity outside a compact set.
\eitem
\end{lemma}
Thus, the set $Y_{\eps'} = \cup_i Y_{\eps'}^i $ is a tubular neighborhood of $image(\tau)$ on which we have nice local coordinate
systems, see Figure \ref{fig:tube}. This will allow us to locally {\it bend} the curve $\tau$. From now on, $\phi^i$ always denotes
the extended diffeomorphism from $\bR^2$ to $\bR^2$.

Next, choose a $C^\infty$ function $w_2: \bR \mapsto [0,1]$ supported in $[-1,1]$ satisfying
\beq \label{eq:defi_w2}
|\hat{w}_2(\omega)| \le c \cdot e^{-|\omega|}, \quad \omega \in \bR
\eeq
and
\[
      w_2(\tfrac{t}{\rho}) + w_2(\tfrac{t-1}{\rho}) = 1\: (-\tfrac12 \leq t \leq 0)
      \quad \mbox{and} \quad
      w_2(\tfrac{t}{\rho}) + w_2(\tfrac{t+1}{\rho}) = 1\:  (0 \leq t \leq \tfrac12) .
\]
Define a smooth partition of unity of $[0, 1]$ using $w_2$ by
\[
   w_{2,i}(t/\rho) = w_2((t-t_i)/\rho) , \qquad 1 \leq i \le m,
\]
and accordingly the distributions
\[
     \cC^i = \int_{t_{i-1}}^{t_{i+1}}  w_{2,i} (t/\rho) \delta_{\tau(t)} dt;
\]
the partition of unity property giving $ \sum_i \cC^i = \cC$.


Now consider the action of $\phi^i$ on the distribution $f$
\[
   (\phi^i)^\star f = f \circ \phi^i .
\]
This action induces a linear transformation on the space of curvelet coefficients. With $\alpha(f)$ the
curvelet coefficients of $f$ and $\beta(f)$ the curvelet coefficients of $(\phi^i)^\star f$, we obtain
a linear operator
\[
    M_{\phi^i} (\alpha(f)) = \beta(f) .
\]

It is by now well-known that diffeomorphisms preserve sparsity of frame coefficients when the frame
is based on parabolic scaling (as with curvelets and shearlets), e.g., by  \cite{Smi98} (see also
\cite[Theorem 6.1, page 219]{CD05b}), for any $0 < p \le 1$,
\[
    \| M_\phi \|_{Op,p} := \max \left\{ \sup_\eta  \| (\langle \gamma_\eta , \phi^\star  \gamma_{\eta'} \rangle )_{\eta} \|_p,
      \sup_{\eta'}  \| (\langle \gamma_\eta , \phi^\star  \gamma_{\eta'} \rangle )_{\eta'}  \|_p \right\} < \infty .
\]

%
%

After having carefully bended the curve pieces, we can also reserve the process and glue them together. Choosing, e.g.,
$\beta_j =  ( \langle \gamma_\eta ,  \cC_j \rangle )_\eta$, from the decomposition $\cC_j = \sum_{i=1}^m \cC^i_j$ we have
\[
   \beta_j =   \sum_{i=1}^m M_{\phi^i} \alpha_j .
\]
This decomposition allows us to relate sparsity of coefficients of the linear singularity to
those of the curvilinear singularity:
\begin{eqnarray}
\label{eq:boundcurvecoeff}
     \| \beta_j \|_p & \leq& m^{1/p} \cdot \left( \max_i    \|M_{\phi^i}\|_{Op,p} \right) \cdot  \| \alpha_j  \|_p  .
\end{eqnarray}

Finally,
we define the very special distribution $w\cL$ we will consider, supported on a line segment
$\{0\} \times [-\rho,\rho]$ by
\[
w\cL = w_2(x_2/\rho) \cdot \delta_0(x_1).
\]
Then we can write
\[
\widehat{w\cL} = \hat{w} \star \hat{\cL},
\]
where
\[
\hat{w} = \hat{w}_2(\rho \xi_2) \cdot \rho \cdot \delta_0(\xi_1)
\quad \mbox{and} \quad
\hat{\cL} = \delta_0(\xi_2).
\]
Thus the action of $w\cL$ on a continuous function $f$ is given by
\[
2\pi \langle w\cL , f \rangle  = \langle \cL , \hat{w} \star \hat{f} \rangle
= \int (\hat{w} \star \hat{f})(\xi_1,0) d\xi_1.
\]
Conceptually, $w\cL$ is a straight curve fragment, which the approach taken in \cite{DK08a} reduced the
analysis of $\cC$ to.

Concluding this approach enables us to consider curvelet coefficients of a linear singularity instead
of a curvilinear singularity with a linear operator mapping one coefficient set onto the other.



\subsection{Estimates for Curvelet Coefficients}

We start by estimating the interaction of curvelet atoms and the interaction of a curvelet atom with a
wavelet atom. These results are proved in \cite{DK08a}.

\begin{lemma}[\cite{DK08a}]
\label{lemm:gamma_gamma_estimate}
For each $N=1,2,\dots$,
there is a constant $c_N$ so that
\[
  |\langle \gamma_{a,b,\theta} , \gamma_{a_0,b_0,\theta_0} \rangle |
  \leq c_N \cdot 1_{\{|\log_2(a/a_0)| < 3\}} \cdot 1_{\{ |\theta - \theta_0| < 10 \sqrt{a_0}  \}} \cdot
    \langle |b - b_0|_{a_0,\theta_0} \rangle^{-N} .
\]
\end{lemma}

\begin{lemma}[\cite{DK08a}]
\label{lemm:gamma_psi_estimate}
For each $N=1,2,\dots$, there is  a constant $c_N$ so that
\[
  |\langle \gamma_{a,b,\theta} , \psi_{a_0,b_0} \rangle |
  \leq c_N \cdot a^{1/4} \cdot 1_{\{|\log_2(a/a_0)| < 3\}}  \cdot
    \langle |b - b_0|_{a,\theta} \rangle^{-N} .
\]
\end{lemma}

We now first analyze curvelet coefficients of the point singularity $\cP$. The technical proof will be
given in Subsection \ref{subsec:proofs_geometrycurvelets}.

\begin{lemma}
\label{lemma:pointcurvelet}
For each $N = 1,2, \dots$, there is a constant $c_N$ so that
\[
   \absip{\cP_j}{\gamma_{a_{j},b,\theta}} \le c_N \cdot a_j^{-1/2} \cdot \langle |D_{1/a_j}b| \rangle^{-N},
   \qquad \forall j \in \bZ\; \forall b, \theta.
\]
\end{lemma}

Next we state two lemmata from \cite{DK08a}, which provide
estimates for the curvelet coefficients of our linear singularity by first considering curvelets, which
are almost aligned with the singularity, and secondly considering the remaining ones.

\begin{lemma}[\cite{DK08a}] \label{lemma:linecurvelet_estimate_close0}
Suppose that $\theta \in [0,\sqrt{a}]$, and set
\[
\tau := \cos \theta \sin \theta (a^{-1}-a^{-2}),\quad d_1^2 = b_1^2(\sigma_2^2-\sigma_1^{-2}\tau),
\]
and
\[
d_2^2
= \left\{ \begin{array}{ccl}
\min\{((\rho-b_2)\sigma_1+\sigma_1^{-1}b_1\tau)^2,((-\rho-b_2)\sigma_1+\sigma_1^{-1}b_1\tau)^2\} & : & b_2-\sigma_1^{-2}b_1\tau \not\in[-\rho,\rho],\\
0 & : & b_2-\sigma_1^{-2}b_1\tau \in[-\rho,\rho],
\end{array} \right.
\]
where
\[
\sigma_1 = (a^{-2} \sin^2 \theta + a^{-1} \cos^2 \theta)^{1/2} \quad \mbox{and} \quad
\sigma_2 = (a^{-1} \sin^2 \theta + a^{-2} \cos^2 \theta)^{1/2}.
\]
Then, for $N=1, 2, \ldots$,
\[
|\langle w\cL ,  \gamma_{a,b,\theta}  \rangle |
\le c_N \cdot a^{-3/4} \cdot \sigma_1^{-1} \cdot \langle d_1 \rangle^{-1}
\cdot \langle|(d_1,\sigma_1 d_2)|\rangle^{2-N}.
\]
In particular, if $\theta=0$,
\[
|\langle w\cL ,  \gamma_{a,b,\theta}  \rangle | \leq
c_N \cdot a^{-1/4} \cdot \langle |b_1/a|\rangle^{-1} \cdot \langle a^{-1} [|b_1|^2 + \min\{(b_2-\rho)^2,(b_2+\rho)^2\}]^{1/2} \rangle^{2-N},
\]
and, if $\theta=0$ and $b_2 \in [-\rho,\rho]$,
\[
|\langle w\cL ,  \gamma_{a,b,\theta}  \rangle | \leq c_N \cdot a^{-1/4} \cdot \langle |b_1/a|\rangle^{1-N}.
\]
\end{lemma}

\begin{lemma}[\cite{DK08a}] \label{lemma:linecurvelet_estimate_other}
Suppose that $\theta \in (\sqrt{a},\pi)$. Then, for $N=1, 2, \ldots$,
\begin{eqnarray*}
\absip{w\cL}{\gamma_{a,b,\theta}}
& \le & c_{L,M} \cdot a^{-1/4} \cdot |\cos \theta| \cdot e^{-\rho \frac{|\sin \theta|}{2a}} \cdot \langle |b_1| \rangle^{-L} \cdot
(a^{1/2} |\sin \theta| + a |\cos \theta|)^L\\
& &  \cdot \langle |b_2| \rangle^{-M} \cdot (\rho + a^{1/2} |\cos \theta| + a |\sin \theta|)^M.
\end{eqnarray*}
\end{lemma}


\subsection{Relation of $\cT_{2,j}^{PS}$ to Significant Coefficients from $\ell_1$ Minimization}
\label{subsec:applyDK09}

Comparing the set of significant coefficients $\cT_{2,j}$ we derive from thresholding with the set
of significant coefficients associated with $\ell_1$ minimization studied in \cite{DK08a} will be
quite beneficial, since it will later on allow us to exploit some of the results from this paper.

To start, we briefly review the definitions and choices made for the significant curvelet coefficients
associated with $\ell_1$ minimization. Recalling the definition of the straight curve fragment $w\cL$
from Section \ref{subsec:bending}, we first define a neighborhood of $WF(w\cL)$ by
\beq \label{eq:defiN}
\cN^{PS}(a,c,\eps') = \{b \in \bR^2 : d_2(b,\{0\} \times [-2\rho,2\rho]) \le c \cdot D_2(a,\eps')\} \times [0,\sqrt{a}],
\eeq
where $c > 0$ is some constant and
\[
D_2(a,\eps') = a^{(1-\eps')} \qquad \mbox{for some } \eps' > 0.
\]
Then the set of significant curvelet coefficients for $w\cL$ was in \cite{DK08a} chosen as
\[
\t{\cS}_{j}(c,\eps') = \{ (j,k,\ell) \in \bigcup_{j'=j-1}^{j+1} \Delta_{j'} : (b_{j,k,\ell},\theta_{j,\ell}) \in \cN^{PS}(a_j,c,\eps')\}.
\]

Let us now first consider the set $\tilde{\cT}_{2,j}$ defined by
\[
\tilde{\cT}_{2,j} = \{ \eta : |\ip{w\cL_j}{\gamma_\eta}| \ge 2^{j(1/4-\eps)}\},
\]
which is related to $\t{\cS}_{j}$ in the following way:

\begin{proposition}
\label{prop:ResultsDK09andHere_wL}
There exist $c_1, c_2 > 0$ and $\eps'_1, \eps'_2 \in (0,\eps)$ such that
\[
\t{\cS}_{j}(c_1,\eps'_1) \subseteq \tilde{\cT}_{2,j} \subseteq \t{\cS}_{j}(c_2,\eps'_2).
\]
\end{proposition}

\noindent {\bf Proof.}
We first prove $\tilde{\cT}_{2,j} \subseteq \t{\cS}_{j}(c_2,\eps'_2)$.
Using the estimate in Lemma \ref{lemma:linecurvelet_estimate_close0} as a model, we obtain the following:
For all $(b_{j,k,\ell},\theta_{j,\ell})$ with $(j,k,\ell) \in \tilde{\cT}_{2,j}$ and $N=1, 2, \ldots$,
we have
\[
 \sigma_1^{\frac{1}{N-2}} \cdot \langle d_1 \rangle^{\frac{1}{N-2}} \cdot \langle|(d_1,\sigma_1 d_2)|\rangle
 \le c_N \cdot 2^{j\frac{1/2+\eps}{N-2}}.
\]
Now Lemma \ref{lemma:linecurvelet_estimate_other} implies that WLOG we only need to consider the case $\theta \in [0,\sqrt{a}]$
due to the rapidly decaying exponential factor. To obtain an estimate for $b$, also WLOG we can assume
that $\theta=0$, which implies
\[
a^{-\frac{1}{2(N-2)}} \cdot  \langle |b_1/a| \rangle^{\frac{1}{N-2}}
\cdot \langle a^{-1} [|b_1|^2 + \min_\pm (b_2\pm\rho)^2]^{1/2} \rangle
 \le c_N \cdot 2^{j\frac{1/2+\eps}{N-2}},
\]
which is equivalent to
\beq \label{eq:split0}
\langle |b_1/a| \rangle^{\frac{1}{N-2}} \cdot \langle a^{-1} [|b_1|^2 + \min_\pm (b_2\pm\rho)^2]^{1/2} \rangle
 \le c_N \cdot 2^{j\frac{\eps}{N-2}}.
\eeq
Since both factors are larger than $1$, we can split this inequality into
\beq \label{eq:split1}
\langle |b_1/a| \rangle^{\frac{1}{N-2}} \le c_N \cdot 2^{j\frac{\eps}{N-2}}
\eeq
and
\beq \label{eq:split2}
 \langle a^{-1} [|b_1|^2 + \min_\pm (b_2\pm\rho)^2]^{1/2} \rangle
 \le c_N \cdot 2^{j\frac{\eps}{N-2}}.
\eeq
From \eqref{eq:split1}, we conclude that
\beq \label{eq:estimateb1}
|b_1| \le  c_N \cdot 2^{-j(1-\eps)},
\eeq
and from \eqref{eq:split2}, we conclude that
\beq \label{eq:estimateb2}
\min_\pm |b_2\pm\rho| \le c_N \cdot 2^{-j(1-\frac{\eps}{N-2})}.
\eeq
Thus, for $c_2$ and $\eps'_2 \in (0,\eps)$ appropriately chosen,
\[
\tilde{\cT}_{2,j} \subseteq \t{\cS}_{j}(c_2,\eps'_2).
\]

The converse inclusion can be derived by substituting \eqref{eq:estimateb1} and \eqref{eq:estimateb2} into
\eqref{eq:split0}. This proves the lemma.
\qed

However, we wish to remind the reader that it is the set of significant curvelet coefficients of the
curvilinear singularity we aim to analyze. For this reason, in the approach presented in \cite{DK08a},
the aforementioned linear operator was exploited to obtain the set of significant curvelet
coefficients of $\cC$ based on the chosen set $\t{\cS}_{j}(c,\eps')$ for $w\cL$. For this,
let $M_{F_j} = (\ip{\gamma_{\eta}}{F_j \star \gamma_{\eta'}})_{\eta,\eta'}$ be the filtering matrix
associated with the filter $F_j$. The `correct' linear operator to consider is defined by the matrix
\[
M_j^i = M_{F_j} \cdot M_{(\phi^i)^{-1}},
\]
and the entries of this matrix will be denoted by $M_j^i(\eta,\eta')$.
Further, we let $t_{\eta',n}$ denote the amplitude of the $n$'th largest element of the $\eta'$'th column.
Now setting
\[
  \cS_{j}^i(c,\eps') = \{ \eta : \eta' \in \t{\cS}_{j}(c,\eps') \mbox{ and }  |M_j^i(\eta,\eta')| > t_{\eta',2^{j\eps'}} \},
\]
the overall cluster set of significant curvelet coefficients of $\cC$ is
\[
   \cS_{2,j}(c,\eps') = \bigcup_i  \cS_{j}^i(c,\eps').
\]

Highly technical and tedious computations (compare \cite[Sec. 7]{DK08a}) -- which we decided to not repeat
here due to their non-intuitive nature -- then imply the following result by using Proposition \ref{prop:ResultsDK09andHere_wL}.

\begin{proposition}
\label{prop:ResultsDK09andHere_C}
There exist $c_1, c_2 > 0$ and $\eps'_1, \eps'_2 \in (0,\eps)$ such that
\[
\cS_{2,j}(c_1,\eps'_1) \subseteq \{ \eta : |\ip{\cC_j}{\gamma_\eta}| \ge 2^{j(1/4-\eps)}\} \subseteq \cS_{2,j}(c_2,\eps'_2).
\]
\end{proposition}

This observation ensures that results from \cite{DK08a} concerning the set of significant curvelet coefficients
are transferable to the situation under consideration in this paper; pleasing news which we intend to take advantage
of.



\subsection{Geometry of $\cT_{2,j}^{PS}$}

Our next goal is to show that instead of considering the set $\cT_{2,j}$ which depends on the residual $\cR_j$
-- typically difficile to handle -- we might consider the `easier-to-handle' set
\[
\{ \eta : |\ip{\cC_j}{\gamma_\eta}| \ge 2^{j(1/4-\eps')}\}
\]
with some control on $\eps'$. This requires a careful analysis of the behavior of the coefficients $\ip{\cR_j}{\gamma_\eta}$,
which are of the following form:

\begin{lemma}
\label{lemma:formRj}
We have
\[
\ip{\cR_j}{\gamma_\eta}
= \ip{\cC_j}{\gamma_\eta} - \sum_{\lambda \in \cT_{1,j}} \ip{\cC_j}{\psi_\lambda} \ip{\psi_\lambda}{\gamma_\eta}
+ \sum_{\lambda \in \cT_{1,j}^c} \ip{\cP_j}{\psi_\lambda} \ip{\psi_\lambda}{\gamma_\eta}.
\]
\end{lemma}

\noindent
{\bf Proof.}
We compute
\[
\ip{\cR_j}{\gamma_\eta}
= \ip{\cP_j}{\gamma_\eta} + \ip{\cC_j}{\gamma_\eta} - \ip{\sum_{\lambda \in \cT_{1,j}}\ip{\cP_j}{\psi_\lambda}\psi_\lambda}{\gamma_\eta}
- \sum_{\lambda \in \cT_{1,j}} \ip{\cC_j}{\psi_\lambda}\ip{\psi_\lambda}{\gamma_\eta}.
\]
Using the fact that $(\psi_\lambda)_\lambda$ is a tight frame, we conclude that
\[
\ip{\cP_j}{\gamma_\eta} - \ip{\sum_{\lambda \in \cT_{1,j}}\ip{\cP_j}{\psi_\lambda}\psi_\lambda}{\gamma_\eta}
= \sum_{\lambda \in \cT_{1,j}^c} \ip{\cP_j}{\psi_\lambda} \ip{\psi_\lambda}{\gamma_\eta},
\]
and the lemma is proved. \qed

The threshold was chosen precisely so that $\ip{\cR_j}{\gamma_\eta} \approx \ip{\cC_j}{\gamma_\eta}$ for all $\eta$
asymptotically, i.e., that the two residuals in Lemma \ref{lemma:formRj} become asymptotically negligible. A quantitative
statement of this consideration is

\begin{proposition}
\label{prop:ResidualDecay}
For any $\delta > 0$,
\[
\Bigg|\sum_{\lambda \in \cT_{1,j}} \ip{\cC_j}{\psi_\lambda} \ip{\psi_\lambda}{\gamma_\eta}
+ \sum_{\lambda \in \cT_{1,j}^c} \ip{\cP_j}{\psi_\lambda} \ip{\psi_\lambda}{\gamma_\eta}\Bigg|
= O(2^{-j(1/4-\delta)}), \qquad j \to \infty.
\]
In particular, we have
\[
\{ \eta : |\ip{\cC_j}{\gamma_\eta}| \ge 2^{j(1/4-(\eps-\delta))}\}
\subseteq \{ \eta : |\ip{\cR_j}{\gamma_\eta}| \ge 2^{j(1/4-\eps)}\}
\subseteq \{ \eta : |\ip{\cC_j}{\gamma_\eta}| \ge 2^{j(1/4-(\eps+\delta))}\}.
\]
\end{proposition}

\noindent
{\bf Proof.}
Let $\delta > 0$ be arbitrary. For proving the first claim, we consider both terms on the LHS separately.
By Lemmata \ref{lemma:curvewavelet} and \ref{lemm:gamma_psi_estimate},
\[
\Bigg|\sum_{\lambda \in \cT_{1,j}} \ip{\cC_j}{\psi_\lambda} \ip{\psi_\lambda}{\gamma_\eta}\Bigg|
\le c \cdot \sum_{\lambda \in \cT_{1,j}} \absip{\psi_\lambda}{\gamma_\eta}\\
\le c \cdot |\cT_{1,j}| \cdot 2^{-j/4}.
\]
Since Lemma \ref{lemma:T1} implies that
\[
|\cT_{1,j}| \le c_N \cdot 2^{j\frac{1-2\eps}{N}},
\]
we obtain
\[
\Bigg|\sum_{\lambda \in \cT_{1,j}} \ip{\cC_j}{\psi_\lambda} \ip{\psi_\lambda}{\gamma_\eta}\Bigg|
\le c_N \cdot 2^{j\frac{1-2\eps}{N}-j/4}.
\]
For $N$ large enough,
\beq \label{eq:R1help1}
\Bigg|\sum_{\lambda \in \cT_{1,j}} \ip{\cC_j}{\psi_\lambda} \ip{\psi_\lambda}{\gamma_\eta}\Bigg| \le c \cdot 2^{-j(1/4-\delta)}, \qquad j \to \infty.
\eeq

Secondly, by Lemmata \ref{lemma:pointwavelet} and \ref{lemma:T1c},
\beq \label{eq:rest1}
\sum_{\lambda \in \cT_{1,j}^c} \absip{\cP_j}{\psi_\lambda}
\le c_N \cdot \sum_{|k| > c_N \cdot 2^{j\frac{1-2\eps}{2N}}} 2^{j/2} \cdot \langle |k| \rangle^{-N}.
\eeq
For $N$ large enough, we have
\beq \label{eq:rest2}
\int_{\{x : |x| > c_N \cdot 2^{j\frac{1-2\eps}{2N}}\}} \langle |x| \rangle^{-N} dx_2 dx_2
\le c_N \cdot 2^{j(1-N)\frac{1-2\eps}{N}} \le c_N \cdot 2^{-j(1/2-\delta)}.
\eeq
By \eqref{eq:rest1} and  \eqref{eq:rest2}, also exploiting Lemma \ref{lemm:gamma_psi_estimate},
\beq \label{eq:R1help2}
\Bigg|\sum_{\lambda \in \cT_{1,j}^c} \ip{\cP_j}{\psi_\lambda} \ip{\psi_\lambda}{\gamma_\eta}\Bigg|
\le c \cdot 2^{j\delta} \cdot 2^{-j/4} \le c \cdot 2^{-j(1/4-\delta)}, \qquad j \to \infty.
\eeq
Now the first claim follows from \eqref{eq:R1help1} and \eqref{eq:R1help2}.

The `in particular'-part can now be derived as a consequence of the first claim by using Lemma \ref{lemma:formRj}.
\qed

Lemma \ref{lemm:WFPinT1} already proved that $WF(\cP) \subseteq \cT_{1,j}^{PS}$. Our last result in this
subsection shows that a similar result holds true for the wavefront set of $\cC$ and the thresholding set
$\cT_{2,j}^{PS}$. These two results will be one main ingredient for proving the separation of wavefront
sets through Single-Pass Alternating Thresholding stated in Theorem \ref{theo:thresholding3}.

\begin{lemma} \label{lemm:WFCinT2}
For $j$ sufficiently large,
\[
\absip{w\cL_j}{\gamma_{j,(0,k_2),0}} \ge 2^{j(1/4-\eps)}\qquad \forall\, k_2 \in 2^j [-\rho,\rho].
\]
Hence, in particular,
\[
WF(\cC) \cap \{(b_{j,k,\ell},\theta_{j,\ell}) : (j,k,\ell) \in \Delta_j\} \subseteq \cT_{2,j}^{PS}.
\]
\end{lemma}

\noindent
{\bf Proof.}
By Parseval,
\begin{eqnarray*}
\absip{w\cL_j}{\gamma_{j,(0,k_2),0}}
& = & 2\pi \absip{\hat{w\cL}_j}{\hat{\gamma}_{j,(0,k_2),0}}\\
& = & 2\pi \cdot 2^{-3j/4} \cdot \int \rho \cdot \hat{w}_2(-\rho \xi_2) W(|\xi|/2^j) V(2^{j/2}\omega) e^{i (k_2/2^{j/2}) \xi_2} d\xi.
\end{eqnarray*}
Apply the change of variables $\zeta = (\xi_1/2^j,\xi_2)$ and $d\zeta = 2^{j} d\xi$,
\beq \label{eq:easierestimate1}
\absip{w\cL_j}{\gamma_{j,(0,k_2),0}}
= 2\pi \rho \cdot 2^{j/4} \cdot \int \hat{w}_2(-\rho \zeta_2) \left[\int W(|\zeta^{(j)}|)
V(2^{j/2}\omega(\zeta)) d\zeta_1\right] e^{i (k_2/2^{j/2}) \zeta_2} d\zeta_2.
\eeq
where $\zeta^{(j)} = (\zeta_1,\zeta_2/2^j)$ and $\omega(\zeta)$ denotes the angular component
of the polar coordinates of $\zeta$. As $j \to \infty$, the integration area is asymptotically
(as $j \to \infty$) of the form
\[
\Xi = ([-2,-1/2] \cup [1/2,2]) \times [-2^{j/2},2^{j/2}].
\]
Letting $L_\delta =2^{j(1/2-\delta)}$, the choice of $W$ and $V$ implies that the dependence of
\[
[-L_\delta,L_\delta] \ni \zeta_2 \mapsto \int W(|\zeta^{(j)}|) V(2^{j/2}\omega(\zeta)) d\zeta_1
\]
on $j$ is asymptotically negligible, and that its absolute value is uniformly bounded from below. Thus,
by \eqref{eq:easierestimate1} and taking the rapid decay condition \eqref{eq:defi_w2} on $\hat{w}_2$ into account,
for some $c > 0$,
\beq \label{eq:easierestimate2}
\absip{w\cL_j}{\gamma_{j,(0,k_2),0}}
\ge c \cdot 2^{j/4} \cdot \int \hat{w}_2(-\rho \zeta_2) e^{i (k_2/2^{j/2}) \zeta_2} d\zeta_2.
\eeq
Finally, again by \eqref{eq:defi_w2}, we can conclude that there exists some $c' > 0$ such
that
\beq \label{eq:easierestimate3}
\int \hat{w}_2(-\rho \zeta_2) e^{i (k_2/2^{j/2}) \zeta_2} d\zeta_2 \ge c', \quad \forall k_2 \in 2^j [-\rho,\rho].
\eeq
Combining \eqref{eq:easierestimate2} and \eqref{eq:easierestimate3}, for sufficiently large $j$,
\beq \label{eq:finalclaim}
\absip{w\cL_j}{\gamma_{j,(0,k_2),0}} \ge 2^{j(1/4-\eps)}, \quad \forall k_2 \in 2^j [-\rho,\rho]
\eeq
which was claimed.

For the `in particular'-part, we first observe that due to Proposition \ref{prop:ResidualDecay}, WLOG we can consider
\[
\{\eta : \absip{\cC_j}{\gamma_\eta} \ge 2^{j(1/4-(\eps+\delta))}\}
\]
for defining $\cT_{2,j}^{PS}$. We then employ the careful bending of the curve as detailed in Section \ref{subsec:bending},
Proposition \ref{prop:ResultsDK09andHere_C}, \cite[Lem. 7.8]{DK08a}, and Proposition \ref{prop:ResultsDK09andHere_wL},
as well as the fact that $WF(w\cL) = \{(0,b) : b \in [-\rho,\rho]\} \times \{0\}$.
This consideration allows us to conclude that the claim follows from \eqref{eq:finalclaim}.
\qed


\section{Asymptotic Separation}
\label{sec:asympsep}


This section is devoted to the analysis around and to the proof of Theorem \ref{theo:thresholding1}.
We first consider an abstract separation setting,
which we will subsequently apply to each filtered version of an image composed of
pointline and curvelike structures.


\subsection{Abstract Separation Estimate for Thresholding}

Suppose we have two tight frames $\Phi_1 = (\phi_{1,i})_i$, $\Phi_2 = (\phi_{2,j})_j$
in a Hilbert space $\cH$, and a signal vector $S \in \cH$. We assume that
all frame vectors are normalized to $c$, say, i.e.,
\[
\norm{\phi_{1,i}}_2 = c \quad \mbox{and} \quad \norm{\phi_{2,j}}_2 = c \quad \mbox{for all } i,j.
\]
We know {\it a priori} that there exists a decomposition
\[
        S = S_1^0 + S_2^0,
\]
where $S_1^0$ is sparse in $\Phi_1$ and $S_2^0$ is sparsely represented in $\Phi_2$.

\begin{figure}[ht]
\centering
\framebox{
\begin{minipage}[h]{5.3in}
\vspace*{0.3cm}
{\sc \underline{Abstract Version of One-Step-Thresholding}}

\vspace*{0.5cm}

{\bf Parameters:}\\[-2ex]
\begin{itemize}
\item Signal $S$.\\[-2ex]
\item Thresholds $t_1$ and $t_2$.
\end{itemize}

\vspace*{0.25cm}

{\bf Algorithm:}\\[-2ex]
\begin{itemize}
\item[1)] {\it Threshold Coefficients with respect to Frame $\Phi_1$:}\\[-2ex]
\begin{itemize}
\item[a)] Compute $c_i = \ip{S}{\phi_{1,i}}$ for all $i$.\\[-2ex]
\item[b)] Apply threshold and set $\cT_{1} = \{ i : |c_i| \ge t_1\}$.\\[-2ex]
\end{itemize}
\item[2)] {\it Reconstruct and Residualize $\Phi_1$-Components:}\\[-2ex]
\begin{itemize}
\item[a)] Compute $S_1^\star = \Phi_1 1_{\cT_1} \Phi_1^T S$.\\[-2ex]
\item[b)] Compute $R = S - S_1^\star = \Phi_1 1_{\cT_1^c} \Phi_1^T S.$\\[-2ex]
\end{itemize}
\item[3)] {\it Threshold Coefficients with respect to Frame $\Phi_2$ of Residual:}\\[-2ex]
\begin{itemize}
\item[a)] Compute $d_j = \ip{R}{\phi_{2,j}}$ for all $j$.\\[-2ex]
\item[b)] Apply threshold and set $\cT_2 = \{ j : |d_j| \ge t_2\}.$\\[-2ex]
\end{itemize}
\item[2)] {\it Reconstruct $\Phi_2$-Components:}\\[-2ex]
\begin{itemize}
\item[a)] Compute $S_2^\star = \Phi_2 1_{\cT_2} \Phi_2^T R.$
\end{itemize}
\end{itemize}

\vspace*{0.25cm}

{\bf Output:}\\[-2ex]
\begin{itemize}
\item Significant thresholding coefficients: $\cT_{1}$ and $\cT_{2}$.\\[-2ex]
\item Approximations to $S_1^0$ and $S_2^0$: $S_1^\star$ and $S_2^\star$.
\end{itemize}
\vspace*{0.01cm}
\end{minipage}
}
\caption{Abstract version of \OneStep Algorithm to decompose $S = S_1^0 + S_2^0$.}
\label{fig:onestepthresholdingAbstract}
\end{figure}

Now we consider an abstract version of \OneStep as explained in
Figure \ref{fig:onestepthresholdingAbstract}.
The following result provides us with an estimate for the $\ell_2$-separation error which
\OneStep causes. Interestingly, both the relative sparsity measure and the
cluster coherence are an essential part of this estimate similar to the analysis
of $\ell_1$ minimization (cf. \cite{DK08a}).

\begin{proposition} \label{prop:mainestimatethresholding}
Suppose that $S$ can be decomposed as $S=S_1^0+S_2^0$.
Let $\cT_1$, $\cT_2$, $S_1^\star$, and $S_2^\star$ be computed via the algorithm \OneStep
(Figure \ref{fig:onestepthresholdingAbstract}),
and assume that each component $S_i^0$ is relatively sparse in $\Phi_i$ with respect to
$\cT_i$, $i=1,2$, respectively, i.e.,
\[
\norm{1_{\cT_1^c} \Phi_1^T S_1^0}_1 + \norm{1_{\cT_2^c} \Phi_2^T S_2^0}_1 \le \delta.
\]
Setting $\mu_c = \mu_c (\cT_2, \Phi_2; \Phi_1)$, we have
\beq \label{eq:errorthresholding}
  \norm{S_1^\star-S_1^0}_2 + \norm{S_2^\star-S_2^0}_2
\le c \cdot \left[ (1+\mu_c) \cdot \norm{1_{\cT_1} \Phi_1^T S_2^0}_1 + (2+\mu_c) \cdot \delta\right].
\eeq
\end{proposition}

This proposition will be proven in Subsection \ref{subsec:proofsasympsep_1}.


\subsection{Application to the Separation of $\cP$ and $\cC$}

We now apply the estimate \eqref{eq:errorthresholding} from Proposition \ref{prop:mainestimatethresholding}
to the following situation: $S$
will be the filtered composition of curves and points $f_j$ with $S_1^0$ being the pointlike part $\cP_j$
and $S_2^0$ the curvelike part $\cC_j$. Our two tight frames of interest, $\Phi_1$ and $\Phi_2$, will be
chosen to be radial wavelets and curvelets, and we notice that these are indeed equal-norm as required
by Proposition \ref{prop:mainestimatethresholding}. Finally the approximation to $\cP_j$ and $\cC_j$
computed by the algorithm \OneStep, i.e., $S_1^\star$ and $S_2^\star$ will be
denoted by $W_j$ and  $C_j$, respectively.

Let $\delta_j$ denote the degree of approximation by thresholded coefficients, i.e., the sum
$\delta_j = \delta_{j,1} + \delta_{j,2}$ of the wavelet approximation error to the point singularity:
\[
   \delta_{j,1} = \sum_{\lambda \in \cT_{1,j}^c} |\langle  \psi_\lambda , \cP_j \rangle| ;
\]
and the curvelet approximation error to the curvilinear singularity:
\[
   \delta_{j,2} = \sum_{\eta \in \cT_{2,j}^c} |\langle \gamma_\eta, \cC_j \rangle|  .
\]
Further let $\mu_c(\cT_{2,j},\Phi_2 ; \Phi_1)$ denote the cluster coherence
\[
(\mu_c)_j = \mu_c(\cT_{2,j},\Phi_2 ; \Phi_1) = \max_\lambda \sum_{\eta \in \cT_{2,j}} \absip{\gamma_\eta}{\psi_\lambda},
\]
the maximal coherence of a wavelet to a cluster of thresholded curvelet coefficients.
We then have

\begin{corollary} \label{coro:thresholdingresult}
Suppose that the sequence of significant thresholding coefficients $(\cT_{1,j})$, and $(\cT_{2,j})$
computed via \OneStep (Figure \ref{fig:onestepthresholdingCP}) has {\bf all}  of the following three properties:
(i) asymptotically negligible cluster coherence:
\[
 (\mu_c)_j =   \mu_c(\cT_{2,j}  , \Phi_2 ; \Phi_1) \goto 0, \qquad j \goto \infty,
\]
(ii) asymptotically negligible cluster approximation error:
\[
  \delta_j =  \delta_{j,1} + \delta_{j,2} = o(\| \cP_j\|_2 + \|\cC_j\|_2) , \qquad j \goto \infty,
\]
(iii) asymptotically negligible energy of the wavelet coefficients of $\cC_j$
on $\cT_{1,j}$:
\[
  \sum_{\lambda \in \cT_{1,j}} \absip{\cC_j}{\psi_\lambda} = o(\| \cP_j\|_2 + \|\cC_j\|_2) , \qquad j \goto \infty.
\]
Then we have asymptotically near-perfect separation:
\[
    \frac{ \|W_j - \cP_j \|_2 +  \|C_j - \cC_j \|_2 }{\|\cP_j\|_2 + \|\cC_j\|_2} \goto 0, \qquad j \goto \infty.
\]
\end{corollary}


\subsection{Proof of Theorem \ref{theo:thresholding1}}
\label{sec:thresholdingtheorem1}

We first recall the following result from \cite{DK08a}.

\begin{lemma}[\cite{DK08a}] \label{lemma:PC}
\[
\| \cP_j\|_2 + \|\cC_j\|_2 = \Omega(2^{j/2}), \qquad j \to \infty.
\]
\end{lemma}

It is now sufficient to show that conditions (i)--(iii) in Corollary \ref{coro:thresholdingresult}
hold true, which is the content of the following four short lemmas. Notice that part (ii) is split into two claims.

\begin{lemma} \label{lemm:mucj}
\[
 (\mu_c)_j = \max_\lambda \sum_{\eta \in \cT_{2,j}} \absip{\gamma_\eta}{\psi_\lambda} \to 0, \qquad j \goto \infty.
\]
\end{lemma}

\noindent {\bf Proof.}
By Proposition \ref{prop:ResidualDecay}, it suffices to prove the
result for
\[
\{ \eta : |\ip{\cC_j}{\gamma_\eta}| \ge 2^{j(1/4-(\eps + \delta))}\}
\]
instead of $\cT_{2,j}$ for $\delta > 0$ arbitrarily small. By Proposition \ref{prop:ResultsDK09andHere_C},
\[
\{ \eta : |\ip{\cC_j}{\gamma_\eta}| \ge 2^{j(1/4-(\eps + \delta))}\} \subset S_{2,j}(c,\eps'),
\]
with $\eps' < \eps + \delta < 1/32$. Now the claim follows from \cite[Lem. 7.7]{DK08a}. \qed

\begin{lemma} \label{lemm:delta1j}
\[
 \delta_{1,j} = \sum_{\lambda \in \cT_{1,j}^c} \absip{\cP_j}{\psi_\lambda} = o(\| \cP_j\|_2 + \|\cC_j\|_2 ) , \qquad j \goto \infty.
\]
\end{lemma}

\noindent {\bf Proof.}
By Lemmata \ref{lemma:pointwavelet} and \ref{lemma:T1c},
\[
\sum_{\lambda \in \Lambda_j \setminus \cT_{1,j}} \absip{\cP_j}{\psi_\lambda}
\le c_N \cdot \sum_{|k| > c_N \cdot 2^{j\frac{1-2\eps}{2N}}} 2^{j/2} \cdot \langle |k| \rangle^{-N}.
\]
For $N$ large enough and $\eps < 1/32$, we have
\[
\int_{\{x : |x| > c_N \cdot 2^{j\frac{1-2\eps}{2N}}\}} \langle |x| \rangle^{-N} dx_2 dx_1
\le c_N \cdot 2^{j(1-N)\frac{1-2\eps}{N}} \le c_N \cdot 2^{-\eps j},
\]
hence, by Lemma \ref{lemma:PC},
\[
\sum_{\lambda \in \Lambda_j \setminus \cT_{1,j}} \absip{\cP_j}{\psi_\lambda}
= o(2^{j/2}) = o(\| \cP_j\|_2 + \|\cC_j\|_2) , \qquad j \goto \infty. \mbox{ \qed}
\]

\begin{lemma} \label{lemm:delta2j}
\[
 \delta_{2,j} = \sum_{\eta \in \Delta_j \setminus \cT_{2,j}} \absip{\cC_j}{\gamma_\eta} = o(\| \cP_j\|_2 + \|\cC_j\|_2) , \qquad j \goto \infty.
\]
\end{lemma}

\noindent {\bf Proof.}
The argumentation is similar to the proof of Lemma \ref{lemm:mucj}, this time using
\cite[Lem. 7.5]{DK08a} instead of \cite[Lem. 7.7]{DK08a}. \qed

\begin{lemma} \label{lemm:partiii}
\[
\sum_{\lambda \in \cT_{1,j}} \absip{\cC_j}{\psi_\lambda} = o(\| \cP_j\|_2 + \|\cC_j\|_2) , \qquad j \goto \infty.
\]
\end{lemma}

\noindent {\bf Proof.}
By Lemmata \ref{lemma:curvewavelet} and \ref{lemma:T1}, and by Lemma \ref{lemma:PC},
\[
\sum_{\lambda \in \cT_{1,j}} \absip{\cC_j}{\psi_\lambda}
\le \sum_{|k| \le c_N \cdot 2^{j\frac{1-2\eps}{2N}}} c
= c_N \cdot c \cdot 2^{j\frac{1-2\eps}{N}}
= o(2^{j/2}) = o(\| \cP_j\|_2 + \|\cC_j\|_2) , \qquad j \goto \infty,
\]
for sufficiently large $N$. \qed

The conditions of Corollary \ref{coro:thresholdingresult} are satisfied, hence Theorem \ref{theo:thresholding1} is proven.


\section{Approximation of the Wavefront Sets}
\label{sec:approxWF}

This section is devoted to proving Theorem \ref{theo:thresholding2}.


\subsection{Proof of Theorem \ref{theo:thresholding2} (i)}
\label{subsec:secondclaim}

By Lemma \ref{lemma:T1}, $(b,\theta) \in \cT_{1,j}^{PS}$ implies that, for each $N = 1,2, \dots$, there is a
constant $c_N$ so that
\[
|b/a_j| \le c_N \cdot 2^{j\frac{1-2\eps}{2N}},
\]
hence
\[
\cT_{1,j}^{PS} \subseteq \{b \in \bR^2 : |b| \le c_N \cdot 2^{j\frac{1-2(N+\eps)}{2N}}\} \times \bP^1.
\]
Thus
\[
d_{PS}(\cT_{1,j}^{PS},WF(\cP))
= d_{PS}(\cT_{1,j}^{PS},\{0\} \times \bP^1)
\le c_N \cdot 2^{j\frac{1-2(N+\eps)}{2N}},
\]
which, for sufficiently large $N$, immediately implies Theorem \ref{theo:thresholding2} (i). \qed


\subsection{Proof of Theorem \ref{theo:thresholding2} (ii)}
\label{subsec:firstclaim}

First we observe that, due to Proposition \ref{prop:ResidualDecay}, WLOG we can consider
\[
\{\eta : \absip{\cC_j}{\gamma_\eta} \ge 2^{j(1/4-(\eps+\delta))}\}
\]
instead of $\cT_{2,j}$ with arbitrarily small $\delta > 0$. From application of Proposition \ref{prop:ResultsDK09andHere_C},
\cite[Lem. 7.8]{DK08a}, and Proposition \ref{prop:ResultsDK09andHere_wL}, it follows that
\[
\cT_{2,j}^{PS} \subseteq \{b \in \bR^2 : |b| \le c_N \cdot 2^{j(\eps'+4(\eps+\delta)-1)}\} \times [0,2^{-j(1/2-(\eps+\delta))}],
\]
for some $c$. A similar conclusion as in the proof of Theorem \ref{theo:thresholding2} (i) then yields
\[
\limsup_{j \to \infty} d_{PS}(\cT_{2,j}^{PS},WF(\cC)) = 0,
\]
which is what was claimed. \qed

\section{Separation of the Wavefront Sets}
\label{sec:sepWF}

This section is devoted to proving Theorem \ref{theo:thresholding3}.


\subsection{A crucial Lemma}

For proving Theorem \ref{theo:thresholding3}, we first state a general lemma on curvelet
synthesis and the associated wavefront set,
which will be later applied to the functions $\sum_j F_j \star W_j$ and $\sum_j F_j \star C_j$.

\begin{lemma}
\label{lemm:curveletsynthesis}
Let $\Omega \subset \bR^2 \times [0,\pi)$ be a compact set in phase space, let
$(T_j)_{j \ge 0}$ be a nested sequence of discrete sets such that $T_j \subseteq \Omega$
for all $j \ge j_0$, and let $(d_{a_j,b,\theta})_{j \ge 0,(b,\theta) \in T_j}$ be a
sequence of complex numbers which satisfies
\beq \label{eq:condition_dj}
|d_{a_j,b,\theta}| = O(a_j^{-m}), \qquad  j \to \infty
\eeq
for some $m > 0$. We further define
\[
g_j = \sum_{(b,\theta) \in T_j} d_{a_j,b,\theta} \, \gamma_{a_j,b,\theta}
\]
and assume that $(g_j)_{j \ge 0}$ is a bounded sequence in the Schwartz space.
Then
\[
WF(\sum_j g_j) \subseteq \Omega.
\]
\end{lemma}

\noindent {\bf Proof.}
Let $(b',\theta') \in \Omega^c$ and consider
\[
\ip{g_j}{\gamma_{a_{j'},b',\theta'}}
= \sum_{(b,\theta) \in T_j} d_{a_j,b,\theta} \ip{\gamma_{a_j,b,\theta}}{\gamma_{a_{j'},b',\theta'}}.
\]
Hence, by Lemma \ref{lemm:gamma_gamma_estimate} and \eqref{eq:condition_dj}, for all $N=1,2,\ldots$,
\[
\absip{g_j}{\gamma_{a_{j'},b',\theta'}}
\le c_N \cdot a_j^{-m} \cdot 1_{\{|\log_2(a_j/a_{j'})| < 3\}} \sum_{(b,\theta) \in T_j}
1_{\{ |\theta - \theta'| < 10 \sqrt{a_{j'}}  \}} \cdot     \langle |b - b'|_{a_{j'},\theta'} \rangle^{-N} .
\]
Thus
\[
\absip{\sum_j g_j}{\gamma_{a_{j'},b',\theta'}}
\le c_N \sum_{j = j'-1}^{j'+1} a_j^{-m} \cdot \sum_{(b,\theta) \in T_j}
1_{\{ |\theta - \theta'| < 10 \sqrt{a_{j'}}  \}} \cdot     \langle |b - b'|_{a_{j'},\theta'} \rangle^{-N} .
\]
Since $(b',\theta') \in \Omega^c$ and $T_j \subseteq \Omega$  for all $(b,\theta) \in T_j$ ($j \ge j_0$),
for any $N = 1,2,...$,we have
\[
\langle |b - b'|_{a_{j'},\theta'} \rangle^{-N} = O(a_{j'}^N), \qquad j' \to \infty.
\]
Since $m$ is fixed, we conclude that
\[
\absip{\sum_j g_j}{\gamma_{a_{j'},b',\theta'}} = O(a_{j'}^N), \qquad j' \to \infty
\]
for any $N = 1,2,...$. Then \cite[Thm. 5.2]{CD05a} implies that $(b',\theta') \not\in WF(\sum_j g_j)$.
\qed

The proof of Theorem \ref{theo:thresholding3} will now be build upon this lemma.


\subsection{Proof of Theorem \ref{theo:thresholding3}}

We start by applying Lemma \ref{lemm:curveletsynthesis} to the situation $g_j=C_j$, $T_j=\cT_{2,j}^{PS}$, and
$d_{a_j,b,\theta}=\ip{f_j}{\gamma_{a_j,b,\theta}}$. Observe that \eqref{eq:condition_dj} is satisfied
by the decay estimates for the curvelet coefficients for $w\cL$, Lemma \ref{lemma:linecurvelet_estimate_close0},
and for $\cP$, Lemma \ref{lemma:pointcurvelet}, and by the bound for the curvelet coefficients of $\cC$,
\eqref{eq:boundcurvecoeff}. $\Omega$ can be chosen as $\cN^{PS}(a,c,\eps')$ with carefully selected $c$
and $\eps'$ due to the considerations in Section \ref{subsec:applyDK09}.
Then Lemma \ref{lemm:curveletsynthesis} together with Theorem \ref{theo:thresholding2} imply
\beq \label{eq:WFCsubset}
WF(\sum_j F_j \star C_j) \subseteq WF(\cC).
\eeq
In a similar way -- by an obvious adaption of Lemma \ref{lemm:curveletsynthesis} -- we can show
\beq \label{eq:WFPsubset}
WF(\sum_j F_j \star W_j) \subseteq WF(\cP).
\eeq

Inclusions \eqref{eq:WFCsubset} and \eqref{eq:WFPsubset} are a significant part of what was claimed, however
a stronger result is true. In order to prove
equality for \eqref{eq:WFPsubset}, it suffices to show that -- since $WF(\cP) = \{0\} \times \mathbf{P}^1$ --  the term
\beq \label{eq:slowdecay1}
\ip{\sum_j F_j \star P_j}{\psi_{j',0}}
= \sum_{j=j'-1}^{j'+1} \sum_{\lambda \in \cT_{1,j}} \ip{f_j}{\psi_{\lambda}}\ip{\psi_{\lambda}}{F_j \star \psi_{j',0}}
\eeq
is of slow decay as $j' \to \infty$, i.e., there exists an $N = 1, 2, ...$ such that this term behaves like $\Omega(a^N)$ as $a \to 0$.
Similarly, for proving equality for \eqref{eq:WFCsubset},  it suffices to show that, for all $(b_{j',k',\ell'},\theta_{j',k'}) \in
WF(\cC)$,  the term
\beq \label{eq:slowdecay2}
\ip{\sum_j F_j \star C_j}{\gamma_{j',k',\ell'}}
= \sum_{j=j'-1}^{j'+1} \sum_{\eta \in \cT_{2,j}} \ip{\cR_j}{\gamma_{\eta}}
\ip{\gamma_{\eta}}{F_j \star \gamma_{j',k',\ell'}},
\eeq
is of slow decay as $j' \to \infty$.
By \cite{CD05a} and the comparable result for wavelets, this then implies that
\[
WF(\cC) \subseteq WF(\sum_j F_j \star C_j) \quad \mbox{and} \quad WF(\cP) \subseteq WF(\sum_j F_j \star W_j),
\]
and, combined with \eqref{eq:WFCsubset} and \eqref{eq:WFPsubset}, the theorem is proved.


We now first show slow decay of the term \eqref{eq:slowdecay1}. For this, we partition the term under consideration
into the following three terms:
\beq \label{eq:slowdecay1_splitting}
\sum_{j=j'-1}^{j'+1} \sum_{\lambda \in \cT_{1,j}} \ip{f_j}{\psi_{\lambda}}\ip{\psi_{\lambda}}{F_j \star \psi_{j',0}}
= T_{11} + T_{12} - T_{13},
\eeq
where
\begin{eqnarray*}
T_{11} & = & \sum_{j=j'-1}^{j'+1} \sum_{\lambda \in \cT_{1,j}} \ip{\cC_j}{\psi_{\lambda}}\ip{\psi_{\lambda}}{F_j \star \psi_{j',0}},\\
T_{12} & = & \ip{\cP_j}{F_j \star \psi_{j',0}},\\
T_{13} & = & \sum_{j=j'-1}^{j'+1} \sum_{\lambda \in \Lambda_j \setminus \cT_{1,j}} \ip{\cP_j}{\psi_{\lambda}}\ip{\psi_{\lambda}}{F_j \star \psi_{j',0}}.
\end{eqnarray*}

We start estimating $T_{11}$. WLOG we can assume that $j=j'$, hence
\[
T_{11} = \sum_{\lambda \in \cT_{1,j}} \ip{\cC_j}{\psi_{\lambda}}\ip{\psi_{\lambda}}{F_j \star \psi_{j,0}}.
\]
By Lemmata \ref{lemm:psi_psi_estimate}, \ref{lemma:curvewavelet}, and \ref{lemma:T1},
\[
|T_{11}| \le c_N \cdot \sum_{|k| \le c_N \cdot 2^{j(1-2\eps)/(2N)}} \langle|k|\rangle^{-N}
\le c_N \cdot \sum_{k = 0}^{\lceil c_N \cdot 2^{j(1-2\eps)/(2N)}\rceil}  k^{1-N}.
\]
Since
\[
\int_0^{\lceil c_N \cdot 2^{j(1-2\eps)/(2N)}\rceil} x^{1-N} dx
\le const,
\]
for sufficiently large $N$, it follows that
\beq \label{eq:slowdecay1T1}
|T_{11}| \le c.
\eeq

Next we estimate $T_{12}$. By Lemma \ref{lemm:WFPinT1}, for sufficiently large $j$,
\beq \label{eq:slowdecay1T2}
T_{12} \ge c \cdot  2^{j/2}.
\eeq

For $T_{13}$, we first observe that WLOG we can assume that $j=j'$, hence
\[
T_{13} = \sum_{\lambda \in \Lambda_j \setminus \cT_{1,j}} \ip{\cP_j}{\psi_{\lambda}}\ip{\psi_{\lambda}}{F_j \star \psi_{j,0}}.
\]
By Lemmata \ref{lemm:psi_psi_estimate}, \ref{lemma:pointwavelet}, and \ref{lemma:T1c},
\[
|T_{13}| \le c_N \cdot 2^{j/2} \cdot \sum_{|k| \ge c_N \cdot 2^{j(1-2\eps)/(2N)}} \langle|k|\rangle^{-2N}
\le c_N \cdot 2^{j/2} \cdot \sum_{k = \lfloor c_N \cdot 2^{j(1-2\eps)/(2N)}\rfloor}^\infty  k^{1-2N}.
\]
Since
\[
\int_{\lfloor c_N \cdot 2^{j(1-2\eps)/(2N)}\rfloor}^\infty x^{1-2N} dx
\le c_N \cdot 2^{-j(1-2\eps)(N-1)/N},
\]
it follows -- by choosing $N=2$ -- that
\beq \label{eq:slowdecay1T3}
|T_{13}| \le c \cdot 2^{-j(1-2\eps)/2}.
\eeq

Applying \eqref{eq:slowdecay1T1}--\eqref{eq:slowdecay1T3} to \eqref{eq:slowdecay1_splitting} implies that
the term $\ip{\sum_j F_j \star P_j}{\psi_{j',0}}$ in \eqref{eq:slowdecay1} behaves like $\Omega(2^{j(1/2-\eps)})$,
hence is of slow decay, which was claimed.

\smallskip

Finally, we prove slow decay of the term \eqref{eq:slowdecay2}. By Propositions \ref{prop:ResidualDecay}
and \ref{prop:ResultsDK09andHere_C}, \cite[Lem. 7.8]{DK08a}, and Proposition \ref{prop:ResultsDK09andHere_wL},
and observing that $WF(w\cL) = \{(0,b) : b \in [-\rho,\rho]\} \times \{0\}$, WLOG we might analyze
\[
\sum_{j=j'-1}^{j'+1} \sum_{\eta \in \tilde{\cT}_{2,j}} \ip{w\cL_j}{\gamma_{\eta}} \ip{\gamma_{\eta}}{F_j \star \gamma_{j',(0,k_2'),0}},
\]
where $k_2' \in 2^j [-\rho,\rho]$. For this, we partition the term under consideration
into the following two terms:
\beq \label{eq:slowdecay2_splitting}
\sum_{j=j'-1}^{j'+1} \sum_{\eta \in \tilde{\cT}_{2,j}} \ip{w\cL_j}{\gamma_{\eta}} \ip{\gamma_{\eta}}{F_j \star \gamma_{j',(0,k_2'),0}}
= T_{21} - T_{22},
\eeq
where
\begin{eqnarray*}
T_{21} & = & \ip{w\cL_j}{F_j \star \gamma_{j',(0,k_2'),0}},\\
T_{22} & = & \sum_{j=j'-1}^{j'+1} \sum_{\eta \in \Delta_j \setminus \tilde{\cT}_{2,j}} \ip{w\cL_j}{\gamma_{\eta}} \ip{\gamma_{\eta}}{F_j \star \gamma_{j',(0,k_2'),0}}.
\end{eqnarray*}

The term $T_{21}$ can be directly estimated by Lemma \ref{lemm:WFCinT2} -- the additional convolution with $F_j$
does not affect the asymptotic behavior -- as
\beq \label{eq:slowdecay2T1}
\absip{w\cL_j}{\gamma_{j,(0,k_2'),0}} \ge 2^{j(1/4-\eps)}\qquad \forall\, k_2' \in 2^j [-\rho,\rho].
\eeq

Next, we analyze $T_{22}$ and first notice that WLOG we can assume that $j=j'$. Thus we are left to estimate
\[
T_{22} = \sum_{(j,k,\ell) \in \Delta_j \setminus \tilde{\cT}_{2,j}} \ip{w\cL_j}{\gamma_{j,k,\ell}} \ip{\gamma_{j,k,\ell}}{F_j \star \gamma_{j,(0,k_2'),0}}
\]
for $k_2' \in 2^j [-\rho,\rho]$. By Lemmata \ref{lemma:linecurvelet_estimate_close0} and \ref{lemm:gamma_gamma_estimate},
and Proposition \ref{prop:ResultsDK09andHere_wL} as well as by the definition of $\cN^{PS}(a,c,\eps')$ in \eqref{eq:defiN},
\begin{eqnarray*}
|T_{22}|
& \le & c_N \cdot 2^{j/4} \cdot \sum_{\{k : \|(k_1/2^{j},k_2/2^{j/2})-(\{0\} \times [-2\rho,2\rho])\|_2 \ge c \cdot 2^{j(\eps-1)}\}} \langle|k-(0,k_2')|\rangle^{-N}\\
& \le & c_N \cdot 2^{j/4} \cdot \sum_{\{k : \|(k_1/2^{j},k_2/2^{j/2})-(\{0\} \times [-2\rho,2\rho])\|_2 \ge c \cdot 2^{j(\eps-1)}\}} |k|^{1-N}.
\end{eqnarray*}

Since, for sufficiently large $N$,
\[
\int_{\{(x_1,x_2) : \|(x_1/2^{j},x_2/2^{j/2})-(\{0\} \times [-2\rho,2\rho])\|_2 \ge c \cdot 2^{j(\eps-1)}\}} |(x_1,x_2)|^{1-N} dx_2 dx_1
\le c \cdot 2^{-2\eps j},
\]
it follows that
\beq \label{eq:slowdecay2T2}
|T_{22}| \le c \cdot 2^{j(1/4-2\eps)}.
\eeq

Applying \eqref{eq:slowdecay2T1} and \eqref{eq:slowdecay2T2} to \eqref{eq:slowdecay2_splitting} implies that
the term $\ip{\sum_j F_j \star C_j}{\gamma_{j',k',\ell'}}$ in \eqref{eq:slowdecay2} behaves like $\Omega(2^{j(1/4-\eps)})$,
hence is of slow decay, which was claimed.
\qed


\section{Proofs}

\subsection{Proof of Results from Section \ref{sec:microlocal}}
\label{subsec:proofs_microlocal}

\subsubsection{Proof of Lemma \ref{lemm:psi_psi_estimate}}

Using Parseval, $\langle \psi_{a,b} , \psi_{a_0,b_0} \rangle  =
     2\pi  \int  \hat{\psi}_{a,b}(\xi) \hat{\psi}_{a_0,b_0}(\xi) d\xi$,
we consider
\[
\int  \hat{\psi}_{a,b}(\xi) \hat{\psi}_{a_0,b_0}(\xi) d\xi
= a^{2}  \int W(ar) \overline{W(a_0r)} e^{-i(b-b_0)\xi} d\xi.
\]
Due to the scaling property, this term is non-zero if and only if $|\log_2(a/a_0)| < 3$. Hence from now on
WLOG we can assume that $a=a_0$. Also, WLOG we may assume that $b_0 = 0$. Applying the change of variables
$\zeta = a \xi$,
\begin{eqnarray*}
\int \hat{\psi}_{a,b}(\xi) \hat{\psi}_{a_0,b_0}(\xi) d\xi
& = & a^{2}  \int |W(ar)|^2 e^{-ib\xi} d\xi\\
& = & \int |W(r)|^2 e^{-i(b/a)\xi} d\xi.
\end{eqnarray*}
Applying integration by parts, for any $k=1, 2, ...$,
\begin{eqnarray*}
\absip{\psi_{a,b}}{\psi_{a_0,b_0}}
& = & 2 \pi \cdot |b/a|^{-k} \left|\int \Delta^k[|W(r)|^2] e^{-i(b/a)\xi} d\xi\right|\\
& \le & 2 \pi \cdot |b/a|^{-k} \int |\Delta^k[|W(r)|^2]| d\xi.
\end{eqnarray*}
Hence
\[
(1+|b/a|^{k}) \cdot \absip{\psi_{a,b}}{\psi_{a_0,b_0}} \le \int |W(r)|^2 d\xi + \int |\Delta^k[|W(r)|^2]| d\xi.
\]
Since the integrand is independent on $a$, and further, for each $k =1, 2, ...$,
\[
\langle|b/a|\rangle^k = (1+|b/a|^2)^\frac{k}{2} \le \frac{k}{2}(1+|b/a|^{k}),
\]
the claim follows.

\subsection{Proofs of Results from Section \ref{sec:geometry_curvelets}}
\label{subsec:proofs_geometrycurvelets}

\subsubsection{Proof of Lemma \ref{lemma:pointcurvelet}}

Using Parseval, $\langle \gamma_{a_j,b,\theta},\cP_j \rangle  =
     2\pi  \int    \hat{\gamma}_{a_j,b,\theta}(\xi) \hat{\cP}_j(\xi) d\xi$,
we consider
\[
   \int   \hat{\gamma}_{a,b,\theta}(\xi) \hat{\cP}_j(\xi)   d\xi \
=    \int   a^{3/4}   W(a r) V((\omega-\theta)/\sqrt{a}) e^{-ib'\xi} \cdot W(a r) \cdot r^{-1/2} d\xi.
\]
Now WLOG we may consider the special case
 $\theta=0$, so that $R_\theta = I$.
We may also assume $b_0 = 0$. Apply the change of variables
$\zeta = D_{a} \xi$ and $d\zeta = a^{3/2} d\xi$,
\begin{eqnarray*}
\int   \hat{\gamma}_{a,b,\theta}(\xi) \hat{\cP}_j(\xi)   d\xi
& = & a^{-3/4} \cdot \int W^2(\norm{\zeta_a}) V(\omega(\zeta_a)/\sqrt{a}) \norm{D_{1/a}\zeta}^{-1/2} e^{-i(D_{1/a}b)'\zeta} d\zeta\\
& = & a^{-1/2} \cdot \int W^2(\norm{\zeta_a}) V(\omega(\zeta_a)/\sqrt{a}) \norm{(a^{-1/2}\zeta_1,\zeta_2)}^{-1/2} e^{-i(D_{1/a}b)'\zeta} d\zeta,
\end{eqnarray*}
where $\zeta_a = (\zeta_1,\sqrt{a}\zeta_2)$ and $\omega(\zeta_a)$ denotes the angular component
of the polar coordinates of $\zeta_a$.
Applying integration by parts, for any $k=1, 2, ...$,
\begin{eqnarray}\nonumber
\lefteqn{\absip{\gamma_{a,b,\theta}}{\cP_j}}\\ \nonumber
& = &2 \pi \cdot a^{-1/2} \cdot |D_{1/a}b|^{-k}
\left|\int \Delta^k [W^2(\norm{\zeta_a}) V(\omega(\zeta_a)/\sqrt{a})\norm{(a^{-1/2}\zeta_1,\zeta_2)}^{-1/2}]e^{-i(D_{1/a}b)'\zeta}
d\zeta\right|\\ \nonumber
& \le & 2\pi \cdot a^{-1/2} \cdot |D_{1/a}b|^{-k}
\int \left|\Delta^k[W^2(\norm{\zeta_a}) V(\omega(\zeta_a)/\sqrt{a})\norm{(a^{-1/2}\zeta_1,\zeta_2)}^{-1/2}]\right| d\zeta.
\end{eqnarray}
Hence
\begin{eqnarray}\nonumber
\lefteqn{(1+|D_{1/a}b|^{k}) \cdot \absip{\gamma_{a,b,\theta}}{\cP_j}}\\ \nonumber
& \le & 2\pi \cdot a^{-1/2} \int \left[ \left|W^2(\norm{\zeta_a})\right|
\left|V(\omega(\zeta_a)/\sqrt{a})\right|\norm{(a^{-1/2}\zeta_1,\zeta_2)}^{-1/2} \right.\\ \label{eq:gammapsi_1}
& & \left. + \left|\Delta^k[W^2(\norm{\zeta_a}) V(\omega(\zeta_a)/\sqrt{a})]\norm{(a^{-1/2}\zeta_1,\zeta_2)}^{-1/2} \right|\right]d\zeta.
\end{eqnarray}
Next we show that, for each $k$, there exists $c_k < \infty$ such that
\begin{eqnarray}\nonumber
\lefteqn{\int \left[ \left|W^2(\norm{\zeta_a})\right|\left|V(\omega(\zeta_a)/\sqrt{a})\right| \norm{(a^{-1/2}\zeta_1,\zeta_2)}^{-1/2}\right.}\\ \label{eq:gammapsi_2}
& & \left. + \left|\Delta^k[W^2(\norm{\zeta_a}) V(\omega(\zeta_a)/\sqrt{a})]\right| \right]d\zeta \le c_k,\quad \forall \;  a>0.
\end{eqnarray}
We have
\[
\frac{\partial}{\partial \zeta_1}W^2(\norm{\zeta_a}) = \frac{\partial}{\partial \zeta_1}W^2(\norm{(\,\cdot\,,\sqrt{a}\zeta_2)})(\zeta_1)
\]
and
\[
\frac{\partial}{\partial \zeta_2}W^2(\norm{\zeta_a}) = \sqrt{a}\cdot\frac{\partial}{\partial \zeta_2}W^2(\norm{(\zeta_1,\sqrt{a}\,\cdot\,)})(\zeta_2).
\]
Hence, by induction, the absolute values of the derivatives of $W^2(\norm{\zeta_a})$ are upper bounded independently of $a$.
Also,
\[
\frac{\partial}{\partial \zeta_1} V(\omega(\zeta_a)/\sqrt{a})
= \frac{\partial}{\partial \zeta_1}V(\omega((\,\cdot\,,\zeta_2)_a)/\sqrt{a})(\zeta_1)\cdot g_1(\zeta,a)
\]
and
\[
\frac{\partial}{\partial \zeta_2} V(\omega(\zeta_a)/\sqrt{a})
= \frac{\partial}{\partial \zeta_2}V(\omega((\zeta_1,\,\cdot\,)_a)/\sqrt{a})(\zeta_1)\cdot g_2(\zeta,a),
\]
and tedious computations show that both $|g_1|, |g_2|$ possess an upper bound independently of $a$. Thus,
by induction, the absolute values of the derivatives
of $V(\omega(\zeta_a)/\sqrt{a})$ are upper bounded independently of $a$. Also, obviously, both
$\frac{\partial}{\partial \zeta_1}\norm{(a^{-1/2}\zeta_1,\zeta_2)}^{-1/2}$ as well as
$\frac{\partial}{\partial \zeta_2}\norm{(a^{-1/2}\zeta_1,\zeta_2)}^{-1/2}$ possess an upper bound independently of $a$.
These observations imply \eqref{eq:gammapsi_2}.

Further, for each $k =1, 2, ...$,
\beq \label{eq:gammapsi_3}
\langle|D_{1/a}b|\rangle^k = (1+|D_{1/a}b|^2)^\frac{k}{2} \le \frac{k}{2}(1+|D_{1/a}b|^{k}).
\eeq

To finish, simply combine \eqref{eq:gammapsi_1}, \eqref{eq:gammapsi_2}, and \eqref{eq:gammapsi_3},
and recall that we chose coordinates so that $\theta = 0$.  Translating back to the case
of general $\theta$ gives the full conclusion. \qed


\subsection{Proofs of Results from Section \ref{sec:asympsep}}

\subsubsection{Proof of Proposition \ref{prop:mainestimatethresholding}}
\label{subsec:proofsasympsep_1}

\noindent {\bf Proof.}
Since $\Phi_1$ is a tight frame,
\begin{eqnarray*}
\norm{S_1^\star-S_1^0}_2
& = & \norm{\Phi_1 1_{\cT_1} \Phi_1^T S-\Phi_1 \Phi_1^T S_1^0}_2\\[1ex]
& = & \norm{\Phi_1 1_{\cT_1} \Phi_1^T S_2^0 - \Phi_1 1_{\cT_1^c} \Phi_1^T S_1^0}_2\\[1ex]
& \le & \norm{\Phi_1 1_{\cT_1} \Phi_1^T S_2^0}_2 + \norm{\Phi_1 1_{\cT_1^c} \Phi_1^T S_1^0}_2.
\end{eqnarray*}
Apply relative sparsity of the subsignal $S_1^0$ and the equal-norm condition
on the tight frame $\Phi_1$,
\beq \label{eq:thresholdingestimate1}
\norm{S_1^\star-S_1^0}_2
\le c \cdot \norm{1_{\cT_1} \Phi_1^T S_2^0}_1 + c \cdot \norm{1_{\cT_1^c} \Phi_1^T S_1^0}_1\\[1ex]
\le c \cdot (\norm{1_{\cT_1} \Phi_1^T S_2^0}_1 + \delta).
\eeq
Next we estimate $\norm{S_2^\star-S_2^0}_2$. We start by using the fact that
$\Phi_2$ is a tight frame and also employ the definition of the residual $R$,
\begin{eqnarray*}
\norm{S_2^\star-S_2^0}_2 \hspace*{-0.23cm}
& = & \hspace*{-0.23cm}  \norm{\Phi_2 1_{\cT_2} \Phi_2^T R-S_2^0}_2\\[1ex]
& = & \hspace*{-0.23cm}  \norm{\Phi_2 1_{\cT_2} \Phi_2^T \Phi_1 1_{\cT_1^c} \Phi_1^T S - \Phi_2 \Phi_2^T \Phi_1 \Phi_1^T S_2^0}_2\\[1ex]
& \le & \hspace*{-0.23cm}  \norm{\Phi_2 1_{\cT_2} \Phi_2^T \Phi_1 1_{\cT_1} \Phi_1^T S_2^0}_2
+ \norm{\Phi_2 1_{\cT_2^c} \Phi_2^T \Phi_1 \Phi_1^T S_2^0}_2
+ \norm{\Phi_2 1_{\cT_2} \Phi_2^T \Phi_1 1_{\cT_1^c} \Phi_1^T S_1^0}_2.
\end{eqnarray*}
Since $\Phi_1$ is a tight frame and the norms of all elements in the tight frame $\Phi_2$ coincide,
we can conclude that
{\allowdisplaybreaks
\begin{eqnarray*}
\lefteqn{\norm{S_2^\star-S_2^0}_2}\\[1ex]
& = & \hspace*{-0.23cm} \norm{\sum_{j \in \cT_2} \sum_{i \in \cT_1} \ip{\phi_{1,i}}{S_2^0} \ip{\phi_{1,i}}{\phi_{2,j}} \phi_{2,j}}_2
+ \norm{\sum_{j \in \cT_2^c} \ip{\phi_{2,j}}{S_2^0} \phi_{2,j}}_2\\
& & + \norm{\sum_{j \in \cT_2} \sum_{i \in \cT_1^c} \ip{\phi_{1,i}}{S_1^0} \ip{\phi_{1,i}}{\phi_{2,j}} \phi_{2,j}}_2\\
& \le & c \cdot \Bigg[
\sum_{i \in \cT_1}\Big(|\ip{\phi_{1,i}}{S_2^0}| \sum_{j \in \cT_2} |\ip{\phi_{1,i}}{\phi_{2,j}}|\Big)
+ \sum_{i \in \cT_1^c}\Big(|\ip{\phi_{1,i}}{S_1^0}|\sum_{j \in \cT_2} |\ip{\phi_{1,i}}{\phi_{2,j}}|\Big)\\
& & +\sum_{j \in \cT_2^c} |\ip{\phi_{2,j}}{S_2^0}| \Bigg].
\end{eqnarray*}
}
Now we have reached the point, where cluster coherence and relative sparsity come into
play. These notions allow us to derive
\begin{eqnarray}\nonumber
\norm{S_2^\star-S_2^0}_2
& \le & c \cdot \Bigg[
\sum_{i \in \cT_1}\Big(|\ip{\phi_{1,i}}{S_2^0}| \cdot \mu_c \Big)
+ \sum_{i \in \cT_1^c}\Big(|\ip{\phi_{1,i}}{S_1^0}| \cdot \mu_c\Big)
+\sum_{j \in \cT_2^c} |\ip{\phi_{2,j}}{S_2^0}|
\Bigg]\\ \label{eq:thresholdingestimate2}
& \le & c \cdot \left[ \mu_c \cdot (\norm{1_{\cT_1} \Phi_1^T S_2^0}_1 + \delta)
+  \delta\right].
\end{eqnarray}
Combining \eqref{eq:thresholdingestimate1} and \eqref{eq:thresholdingestimate2} proves the lemma.
\qed


\end{document}